# Forward Completeness and Applications to Control of Automated Vehicles

**Iasson Karafyllis[*], Dionysis Theodosis[**], and Markos Papageorgiou[**],[***]**

[*]Dept. of Mathematics, National Technical University of Athens,
Zografou Campus, 15780, Athens, Greece,
emails: iasonkar@central.ntua.gr , iasonkaraf@gmail.com

[**]Dynamic Systems and Simulation Laboratory,
Technical University of Crete, Chania, 73100, Greece
(emails: dtheodosis@dssl.tuc.gr , markos@dssl.tuc.gr)

[***]Faculty of Maritime and Transportation,
Ningbo University, Ningbo, China

**Abstract**

Forward complete systems are guaranteed to have solutions that exist globally for all positive time. In this paper, a relaxed Lyapunov-like condition for forward completeness is presented for finite-dimensional systems defined on open sets that does not require boundedness of the Lyapunov-like function along the solutions of the system. The corresponding condition is then exploited for the design of autonomous two-dimensional movement, with focus on lane-free cruise controllers for automated vehicles described by the bicycle kinematic model. The derived feedback laws (cruise controllers) are decentralized and can account for collision avoidance, roads of variable width, on-ramps and off-ramps as well as different desired speed for each vehicle.

## 1. Introduction

In many problems in Physics and Engineering, dynamical systems are subject to various state constraints, and thus, their state space is not a linear space but rather open or closed sets [4]. Inevitably, the control of such systems presents significant difficulties since state constraints need also to be taken into account in the control design. Among the various approaches that tackle systems with state constraints are those of (i) Barrier function-based approaches (see for instance [1], [2], [11], [13], [41]); and (ii) Lyapunov-based approaches (see [5], [12], [15], [16], [17], [33], [34], [39]).

One of the control problems that need to be studied on a set which is not diffeomorphic to a linear space is the design problem of cruise controllers (feedback laws) for automated vehicles. Due to various constraints, such as collision avoidance and speed positivity, the overall system of automated vehicles under consideration is defined on an open set (see [15], [16], [17], [19], [34], [38]). The problem has characteristics that are rarely met in other control problems: besides the fact that the system is defined on an open set (which is not diffeomorphic to $\mathbb{R}^n$), the set of equilibrium points for this system may be not compact. The latter feature is a consequence of the fact that automated vehicles under the effect of decentralized cruise controllers do not interact with each other when the distance between them is large.

The contribution of the present paper is twofold. First, we provide a relaxed Lyapunov-like condition for forward completeness of finite-dimensional systems defined on open sets (Theorem 1). The corresponding Lyapunov-like condition is analogous to the results in [3] for systems defined on $\mathbb{R}^n$ and to the results in [23] for abstract systems defined on normed linear spaces, and does not require boundedness of the Lyapunov-like function along the trajectories of the system (see for instance [18]). Lyapunov stability theory on open sets has a long history and leads to advanced stability notions, see for instance [6], [8], [22], [32].

The second contribution of the paper concerns the design of cruise-controllers (feedback laws) for autonomous two-dimensional movement, with focus on autonomous vehicles driving on lane-free roads (see also [16], [26], [34], [40]). Based on potential functions and additional penalty terms, we construct a Lyapunov-like function which allows for an explicit cruise controller design that guarantees the aforementioned relaxed Lyapunov-like condition for forward completeness (Theorem 2). The Lyapunov-like function acts also as a size function (see [30]) preventing escape of the state from the state-space. The main features of the proposed lane-free cruise controllers are:

- Each vehicle may move within the general (possibly curved) road boundaries or within desired individual corridor boundaries;
- vehicles do not collide with each other;
- the speeds of all vehicles are always positive and remain below a given speed limit;
- each vehicle can have a different longitudinal speed set-point and different size;
- the cruise controllers are fully decentralized, and each vehicle only needs access to the distance from the (curved) boundaries of the road and relative speeds and distances from adjacent vehicles; and
- all the above features are valid globally, i.e., for all physically relevant initial conditions.

If moreover, all vehicles are moving on a road or corridor of constant width and all vehicles have the same speed set-point, then it is shown that all vehicle speeds converge asymptotically to the speed set-point (Theorem 3), thus retaining all the features of the controllers in [15], [16] for the constant-width case.

While the proposed cruise controllers have the same properties of the controllers in [15] and [16] (for the case of constant-width roads), their design is based on a different Lyapunov-like function, which allows for a significant broadening of application conditions. For example, vehicles can now have different speed set-points and move within individual (possibly curved) corridors which allows for the presence of (and driving on) on-ramps, off-ramps, and weaving sections, whereby each vehicle only uses relative speed and relative position measurements from adjacent vehicles to enter or exit the highway and does not require communication with a coordinator or supervisory control algorithms that is common in the literature for lane-based traffic (see for instance [10], [29], [42] and references therein). To our knowledge, cruise controllers that achieve all the above properties simultaneously and can handle most of the possible cases that appear in practice are not available in the literature. In fact, the proposed setting has also the potential to be used for automated driving of vehicles with lane-keeping features (see [2], [7], [27], [36]) or for lane-changing manoeuvres (see for instance [9], [14], [37]), as well as for navigation and control of mobile robots (see [21], [35]), vessels [25], [31] and more.

While the design of the cruise controllers for roads of non-constant width guarantees all above-mentioned features, it does not guarantee a positive lower bound for the speed of all vehicles. Simulations show that in some very special cases, vehicles may asymptotically come to a halt (see also [21] for analogous topological obstructions in robotic systems).

The structure of the paper is as follows. Section 2 is devoted to forward completeness of finite dimensional systems defined on open sets. Section 3 presents the cruise controllers for automated vehicles on roads of variable width. Numerical examples are given in Section 4 illustrating the properties of the proposed cruise controllers. All proofs of the main results are provided in Section 5. Finally, some concluding remarks are given in Section 6.

**Notation.** Throughout this paper, we adopt the following notation.

* $\mathbb{R}_+ := [0, +\infty)$ denotes the set of non-negative real numbers.

* By $|x|$ we denote both the Euclidean norm of a vector $x \in \mathbb{R}^n$ and the absolute value of a scalar $x \in \mathbb{R}$. By $x'$ we denote the transpose of a vector $x \in \mathbb{R}^n$. By $|x|_\infty = \max\{|x_i|, i=1,...,n\}$ we denote the infinity norm of a vector $x = (x_1, x_2, ..., x_n)' \in \mathbb{R}^n$.

* Let $A \subseteq \mathbb{R}^n$ be an open set. By $C^0(A, \Omega)$, we denote the class of continuous functions on $A \subseteq \mathbb{R}^n$, which take values in $\Omega \subseteq \mathbb{R}^m$. By $C^k(A; \Omega)$, where $k \geq 1$ is an integer, we denote the class of functions on $A \subseteq \mathbb{R}^n$ with continuous derivatives of order $k$, which take values in $\Omega \subseteq \mathbb{R}^m$. When $\Omega = \mathbb{R}$ then we write $C^0(A)$ or $C^k(A)$. For a function $V \in C^1(A; \mathbb{R})$, the gradient of $V$ at $x \in A \subseteq \mathbb{R}^n$, denoted by $\nabla V(x)$, is the row vector $\left[\frac{\partial V}{\partial x_1}(x) \cdots \frac{\partial V}{\partial x_n}(x)\right]$.

## 2. A Relaxed Lyapunov-like Condition for Forward Completeness

In this work we consider time-invariant, finite-dimensional systems defined on open sets, i.e., systems of the form

$$\dot{w} = f(w), w \in \Omega \qquad (2.1)$$

where $\Omega \subseteq \mathbb{R}^n$ is a non-empty open set and $f : \Omega \to \mathbb{R}^n$ is a locally Lipschitz mapping. Forward completeness for systems of the form (2.1), i.e., the property that guarantees that for every $w_0 \in \Omega$ the unique solution $w(t) \in \Omega$ of the initial-value problem (2.1) with initial condition $w(0) = w_0$ exists for all $t \geq 0$ has been studied in the literature by means of Lyapunov-like functions (see [18]). Important role in the study of forward completeness plays the notion of the "size function" for system (2.1), i.e., a continuous function $\rho : \Omega \to \mathbb{R}_+$ for which there exists a compact set $A \subset \Omega$ with $\rho(A) = 0$ and $\rho(w) > 0$ for $w \in \Omega \setminus A$ and such that the sublevel sets $\{w \in \Omega : \rho(w) \leq r\}$ are compact for every $r \geq 0$ (see [30]).

However, for systems of the form (2.1) that either do not have equilibrium points or the set of equilibrium points is not bounded, the construction of a radially unbounded Lyapunov-like function for (2.1), i.e., a $C^1$ function $H : \Omega \to \mathbb{R}_+$ that satisfies $\nabla H(w) f(w) \leq 0$ and $R + H(w) \geq \rho(w)$ for all $w \in \Omega$ for a size function $\rho$ and a constant $R \geq 0$ (see [18]), is not easy. This is particularly true for systems of autonomous vehicles which are not restricted to be located within a compact set or more generally for many-body problems in unrestricted space; see next section for more details. To this purpose, we provide the following relaxed Lyapunov-like condition which guarantees forward completeness for systems of the form (2.1).

**Theorem 1:** *Let $\Omega \subseteq \mathbb{R}^n$ be an open set and $f : \Omega \to \mathbb{R}^n$ be a locally Lipschitz vector field. Let $H \in C^1(\Omega; \mathbb{R}_+)$ be a function for which there exist constants $c, \sigma \geq 0$ such that the following inequality holds for all $w \in \Omega$:*

$$\nabla H(w) f(w) \leq cH(w) + \sigma \qquad (2.2)$$

*Furthermore, for each* $r > \inf\{H(w): w \in \Omega\}$ *define the set*

$$S_r = \{w \in \Omega : H(w) \leq r\} \tag{2.3}$$

*and suppose that* $f(S_r)$ *is bounded for every* $r > \inf\{H(w): w \in \Omega\}$. *Furthermore, suppose that for every* $r > \inf\{H(w): w \in \Omega\}$ *and for every* $\rho > \inf\{|w| : w \in S_r\}$ *the set*

$$Q_{r,\rho} = \{w \in \Omega : H(w) \leq r, |w| \leq \rho\} \tag{2.4}$$

*is compact. Then the dynamical system* $\dot{w} = f(w)$ *is forward complete, i.e., for every* $w_0 \in \Omega$ *the unique solution* $w(t) \in \Omega$ *of the initial-value problem* $\dot{w} = f(w)$ *with* $w(0) = w_0$ *is defined for all* $t \geq 0$.

**Discussion of Theorem 1:**
- Notice that when $H \in C^1(\Omega; \mathbb{R}_+)$ satisfies (2.2) and there exist a size function $\rho$ and a constant $R \geq 0$ for which $R + H(w) \geq \rho(w)$ for all $w \in \Omega$ then all conditions of Theorem 1 are automatically satisfied. In particular, in this case the set $S_r$ defined by (2.3) is compact for every $r > \inf\{H(w): w \in \Omega\}$ and consequently (by continuity) the set $f(S_r)$ defined by (2.3) is compact for every $r > \inf\{H(w): w \in \Omega\}$.
- Notice that Theorem 1 does not require the demanding inequality $\nabla H(w) f(w) \leq 0$ for all $w \in \Omega$ but rather the less demanding inequality (2.2). This puts Theorem 1 in the same spirit with the results in [3] for systems defined on $\mathbb{R}^n$ and with the results in [23] for abstract systems defined on normed linear spaces.
- The usefulness of Theorem 1 is shown in the following section for systems of automated vehicles which are not restricted to be located within a compact set. For such systems, it is difficult to construct Lyapunov-like functions or even size functions for the study of forward completeness. On the other hand, for such systems, functions that are constructed by means of "energy" arguments satisfy all requirements of Theorem 1.
- Theorem 1 can be easily extended to the case where $H \in C^0(\Omega; \mathbb{R}_+)$ is not necessarily continuously differentiable but is only locally Lipschitz on $\Omega$.

We next provide a proof of Theorem 1.

**Proof of Theorem 1:** Let $w_0 \in \Omega$ be given. Since $f$ is locally Lipschitz, the initial-value problem $\dot{w} = f(w)$ with initial condition $w(0) = w_0$ has a unique solution $w(t) \in \Omega$ defined for $t \in [0, t_{\max})$, where $t_{\max} \in (0, +\infty]$. Furthermore, if $t_{\max} < +\infty$ then, Theorem 3 on page 91 in [24] shows that for every compact set $K \subset \Omega$ there exists a time $t \in (0, t_{\max})$ for which $w(t) \notin K$. Note that while Theorem 3 in [24] is stated with $f \in C^1(\Omega)$, it can directly be extended with $f$ being locally Lipschitz. We show next by means of a contradiction argument that $t_{\max} = +\infty$.

Without loss of generality, we may assume that (2.2) holds with $c, \sigma > 0$. Suppose that $t_{\max} < +\infty$. We show next that there exists a compact set $K \subset \Omega$ for which $w(t) \in K$ for all $t \in [0, t_{\max})$.

Using (2.2) and the Comparison Principle (see [23]), we obtain that

$$H(w(t)) \leq \exp(ct) H(w(0)) + \frac{\sigma}{c}(\exp(ct) - 1) \text{ for all } t \in [0, t_{\max}) \tag{2.5}$$

which implies that

$$H(w(t)) \leq \exp(ct_{max})H(w_0) + \frac{\sigma}{c}(\exp(ct_{max}) - 1) \tag{2.6}$$

Let $r = \exp(ct_{max})H(w_0) + \frac{\sigma}{c}(\exp(ct_{max}) - 1)$ and notice that definition (2.3) and (2.6) imply that $w(t) \in S_r$ for all $t \in [0, t_{max})$. Since by hypothesis, $f(S_r)$ is bounded, it follows that there exists a constant $M > 0$, such that

$$|\dot{w}(t)| \leq M \quad \text{for all } t \in [0, t_{max}) \tag{2.7}$$

Inequality (2.7) implies that

$$|w(t) - w_0| \leq M t, \text{ for all } t \in [0, t_{max})$$

and consequently (using the triangle inequality)

$$|w(t)| \leq |w_0| + M t_{max} \tag{2.8}$$

Therefore (by virtue of (2.6), (2.8) and definition (2.4)) $w(t) \in K$ for all $t \in [0, t_{max})$, where $K := Q_{r,\rho} \subset \Omega$ with $r = \exp(ct_{max})H(w_0) + \frac{\sigma}{c}(\exp(ct_{max}) - 1)$ and $\rho = |w_0| + M t_{max}$. By assumption, the set $K := Q_{r,\rho} \subset \Omega$ is compact. The proof is complete. ◁

## 3. Autonomous Two-Dimensional Movement of Vehicles in Roads of Variable Width

In this section we illustrate the Lyapunov-like condition of Theorem 1 for forward completeness for the design of lane-free cruise controllers for automated vehicles on roads of variable width. We consider $n$ vehicles described by the bicycle kinematic model, [28]:

$$\begin{aligned}
\dot{x}_i &= v_i \cos(\theta_i) \\
\dot{y}_i &= v_i \sin(\theta_i) \\
\dot{\theta}_i &= \sigma_i^{-1} v_i \tan(\delta_i) \\
\dot{v}_i &= F_i
\end{aligned} \tag{3.1}$$

for $i = 1, \ldots, n$, where $\sigma_i > 0$ is the length of each vehicle (a constant). Here, $(x_i, y_i)$ is the reference point of the $i$-th vehicle and is placed at the midpoint of the rear axle, with $x_i$ being the longitudinal position and $y_i$ being the lateral position of the vehicle in an inertial frame with Cartesian coordinates $(X, Y)$; $v_i > 0$ is the speed of the $i$-th vehicle at the point $(x_i, y_i)$, $\theta_i \in \left(-\frac{\pi}{2}, \frac{\pi}{2}\right)$ is the orientation of the $i$-th vehicle with respect to the $X$ axis, $\delta_i$ is the steering angle of the front wheels relative to the orientation $\theta_i$ of the $i$-th vehicle, and $F_i$ is the acceleration of the $i$-th vehicle. To simplify the subsequent analysis, we use the feedback transformation

$$\delta_i = \arctan\left(\frac{\sigma_i u_i}{v_i}\right), \quad i = 1, \ldots, n \tag{3.2}$$

to obtain the model

$$\dot{x}_i = v_i \cos(\theta_i)$$
$$\dot{y}_i = v_i \sin(\theta_i)$$
$$\dot{\theta}_i = u_i \quad (3.3)$$
$$\dot{v}_i = F_i$$

for $i = 1,...,n$, where $u_i$ and $F_i$, are the inputs of the system.

Let $v_{\max} > 0$ be the road speed limit and $\varphi \in \left(0, \frac{\pi}{2}\right)$ be an angle that satisfies

$$\cos(\varphi) > \max\left\{\frac{\max_{i=1,...,n}\{v_i^*\}}{v_{\max}}, \frac{1}{3}\right\}. \quad (3.4)$$

where $v_i^* \in (0, v_{\max})$ is the speed set-point of each vehicle $i = 1,...,n$. It should be noticed that in contrast to [16], we consider here individual speed set-points $v_i^* \in (0, v_{\max})$ for $i = 1,...,n$. Inequality (3.4) is a technical condition to restrict the maximum orientation angle of the vehicles, $\theta_i \in (-\varphi, \varphi)$.

Suppose that each vehicle $i = 1,...,n$ follows a path which is contained in the set (corridor)

$$Y_i = \{(x, y) : x \in \mathbb{R}, \alpha_i(x) < y < \beta_i(x)\} \quad (3.5)$$

where $\alpha_i : \mathbb{R} \to \mathbb{R}$ and $\beta_i : \mathbb{R} \to \mathbb{R}$ are $C^2$ functions with the following property: there exist constants $r_{\max} > r_{\min} > 0$ and $\gamma_{\max}, \gamma_{\min} \in \mathbb{R}$ with $\gamma_{\max} > \gamma_{\min}$ such that

$$\gamma_{\max} \geq \beta_i(x), \quad \alpha_i(x) \geq \gamma_{\min} \text{ for all } x \in \mathbb{R}, \; i = 1,...,n \quad (3.6)$$

$$r_{\max} \geq \beta_i(x) - \alpha_i(x) \geq r_{\min} > 0 \text{ for all } x \in \mathbb{R}, \; i = 1,...,n \quad (3.7)$$

$$\max\left(\sup_{x \in \mathbb{R}}(|\alpha_i'(x)|), \sup_{x \in \mathbb{R}}(|\beta_i'(x)|)\right) < \tan(\varphi) \text{ for all } x \in \mathbb{R}, \; i = 1,...,n \quad (3.8)$$

$$\sup_{x \in \mathbb{R}}(|\alpha_i''(x)| + |\beta_i''(x)|) < +\infty \text{ for all } x \in \mathbb{R}, \; i = 1,...,n \quad (3.9)$$

Inequality (3.8) describes the maximum rate of change of each path and ensures that the rate of change will not be faster than the maximum angle of orientation of each vehicle (recall (3.4)). Inequalities (3.7) describe the minimum and maximum width of each corridor $(\alpha_i(x), \beta_i(x))$, $i = 1,...,n$. Moreover, the path of each vehicle can be in different sets (corridors) $Y_i = \{(x, y) : x \in \mathbb{R}, \alpha_i(x) < y < \beta_i(x)\}$ for $i = 1,...,n$, and the set $Y_i$ of a vehicle may or may not have common points with the set $Y_j$ with $j \neq i$ of a different vehicle. When all vehicles share the same set $Y_i$, namely, when $\alpha_i(x) = \alpha(x)$ and $\beta_i(x) = \beta(x)$, for all $i = 1,...,n$, then, the curves $\alpha(x)$ and $\beta(x)$ describe the right and left boundary of the road, respectively.

We define the distance between vehicles by

$$d_{i,j} := \sqrt{(x_i - x_j)^2 + p_{i,j}(y_i - y_j)^2}, \text{ for } i, j = 1,...,n \quad (3.10)$$

where $p_{i,j} > 0$ are positive constants that satisfy $p_{i,j} = p_{j,i}$, $i, j = 1...,n$, $j \neq i$. Note that, for $p_{i,j} = 1$, we obtain the standard Euclidean distance, while for $p_{i,j} > 1$, we have an "elliptical" metric, which can approximate more accurately the dimensions of each vehicle. Note that for the case of $n$ vehicles of equal length, the optimal selection of a single $p$ can be found in [15].

Let
$$w = (x_1,...,x_n, y_1,...,y_n, \theta_1,...,\theta_n, v_1,....,v_n)' \in \mathbb{R}^{4n}. \tag{3.11}$$

For each vehicle $i = 1,...,n$ the following constraints should hold $v_i \in (0, v_{max})$, $\theta_i \in (-\varphi, \varphi)$, and $y_i \in (\alpha_i(x_i), \beta(x_i))$, $x_i \in \mathbb{R}$. Define also the set

$$\Omega := \{w \in \mathbb{R}^{4n} : x_i \in \mathbb{R}, y_i \in (\alpha_i(x_i), \beta_i(x_i)), v_i \in (0, v_{max}), \theta_i \in (-\varphi, \varphi), d_{i,j} > L_{i,j}, i, j = 1,...,n, j \neq i\} \tag{3.12}$$

where $L_{i,j}$, $i, j = 1,...,n$, $i \neq j$, are positive constants and represent the minimum distance between a vehicle $i$ and a vehicle $j$, with $L_{i,j} = L_{j,i}$ for $i, j = 1,...,n$, $i \neq j$. Notice that the set $\Omega$ in (3.12) is an open set and constitutes the state-space of our system.

Let $V_{i,j} : (L_{i,j}, +\infty) \to \mathbb{R}_+$ $i, j = 1,...,n$, $j \neq i$ be $C^2$ functions that satisfy the following properties

$$\lim_{d \to L_{i,j}^+} (V_{i,j}(d)) = +\infty \tag{3.13}$$

$$V_{i,j}(d) = 0, \text{ for all } d \geq \lambda \tag{3.14}$$

$$V_{i,j}(d) = V_{j,i}(d), \ i, j = 1,...,n, j \neq i \tag{3.15}$$

where $\lambda$ is a positive constant that satisfies

$$\lambda > \max\{L_{i,j}, i, j = 1,...,n, i \neq j\}. \tag{3.16}$$

We also consider the $C^2$ functions $U_i : (-1,1) \to \mathbb{R}_+$, $i = 1,...,n$ that satisfy

$$U_i(0) = 0 \tag{3.17}$$

$$\lim_{s \to (-1)^+} (U_i(s)) = +\infty, \ \lim_{s \to 1^-} (U_i(s)) = +\infty \tag{3.18}$$

The families of functions $V_{i,j}$ and $U_i$ in (3.13)-(3.15) and (3.17), (3.18), respectively, are potential functions, which have been widely used to avoid collisions between vehicles and road boundary violation (see for instance [16], [38]). Condition (3.15) implies that if a vehicle $i$ exerts a force to vehicle $j$, then vehicle $j$ exerts the opposite force to vehicle $i$. The constant $\lambda$ determines the required real-time information radius around each vehicle. When the distance between vehicles is greater than $\lambda$, then condition (3.14) implies that there is no information exchange between vehicles.

Let $b : \mathbb{R} \to (0,1]$ be a non-increasing $C^2$ function that satisfies

$$b(x) = 1 \text{ for } x \leq \varepsilon \text{ and } b(x)x < M \text{ for } x \geq \varepsilon \tag{3.19}$$

$$\sup_{x \in \mathbb{R}}(|b'(x)| + |b''(x)|) < +\infty \tag{3.20}$$

where $M$ and $\varepsilon$ are positive constants. An example of such a function is

$$b(x) = \begin{cases} 1 & x \leq \varepsilon \\ \exp(-(x-\varepsilon)^3) & x > \varepsilon \end{cases}$$

for some positive constant $\varepsilon > 0$.

Let $R > 0$, and define

$$H(w) = \frac{1}{2}\sum_{i=1}^{n} \frac{(v_i \cos(\theta_i) - f_i(w))^2}{(v_{max} - v_i)v_i} + \frac{R}{2}\sum_{i=1}^{n} \frac{(\sin(\theta_i) - g_i(x_i, y_i)\cos(\theta_i))^2}{\cos(\theta_i) - \cos(\varphi)}$$
$$+ \sum_{i=1}^{n} U_i\left(\frac{2y_i - (\beta_i(x_i) + \alpha_i(x_i))}{\beta_i(x_i) - \alpha_i(x_i)}\right) + \frac{1}{2}\sum_{i=1}^{n}\sum_{j \neq i} V_{i,j}(d_{i,j}) \quad (3.21)$$

where

$$g_i(x, y) := \frac{(\beta_i(x) - y)\alpha_i'(x) + (y - \alpha_i(x))\beta_i'(x)}{\beta_i(x) - \alpha_i(x)}, \text{ for all } x \in \mathbb{R}, \ y \in (\alpha_i(x), \beta_i(x)) \quad (3.22)$$

and $f_i$, $i = 1,...,n$ are $C^1$ functions defined by

$$f_i(w) = v_i^* b(\Phi_i(w) + g_i(x_i, y_i)\Xi_i(w)), \ w \in \Omega, \ i = 1,...,n \quad (3.23)$$

with $\Phi_i(w)$ and $\Xi_i(w)$ for $i = 1,...,n$ given by

$$\Phi_i(w) := \sum_{j \neq i} V_{i,j}'(d_{i,j}) \frac{(x_i - x_j)}{d_{i,j}}, \text{ for } w \in \Omega \quad (3.24)$$

$$\Xi_i(w) = \sum_{j \neq i} p_{i,j} V_{i,j}'(d_{i,j}) \frac{(y_i - y_j)}{d_{i,j}}, \text{ for } w \in \Omega. \quad (3.25)$$

Notice that due to definitions (3.19), (3.23), $f_i(w) \in (0, v_i^*]$ for all $w \in \Omega$.

Let

$$F_i = \frac{1}{q(v_i, \theta_i, f_i(w))}\left(-\mu_1(v_i \cos(\theta_i) - f_i(w)) + \frac{Z_i(w) + v_i \sin(\theta_i)u_i}{v_i(v_{max} - v_i)} - \Phi_i(w) - g_i(x_i, y_i)\Xi_i(w)\right) \quad (3.26)$$

and

$$u_i = \frac{2(\cos(\theta_i) - \cos(\varphi))^2}{Rh_i(x_i, y_i, \theta_i)}\left(-\mu_2 v_i^2(\sin(\theta_i) - g_i(x_i, y_i)\cos(\theta_i))\right.$$
$$\left. + R\frac{v_i \cos(\theta_i)a_i(x_i, y_i, \theta_i)}{\cos(\theta_i) - \cos(\varphi)} - U_i'\left(\frac{2y_i - (\beta_i(x_i) + \alpha_i(x_i))}{\beta_i(x_i) - \alpha_i(x_i)}\right)\frac{2v_i}{\beta_i(x_i) - \alpha_i(x_i)} - v_i \Xi_i(w)\right) \quad (3.27)$$

where

$$Z_i(w) := v_i^* b'(\Phi_i(w) + g_i(x_i, y_i)\Xi_i(w)) \frac{d}{dt}(\Phi_i(w) + g_i(x_i, y_i)\Xi_i(w)), \ i = 1,...,n, \ w \in \Omega \quad (3.28)$$

$$q(v, \theta, r) = \frac{v_{max} v \cos(\theta) + rv_{max} - 2rv}{2(v_{max} - v)^2 v^2}, \text{ for all } v \in (0, v_{max}), \theta \in (-\varphi, \varphi), r \in \mathbb{R} \quad (3.29)$$

$$h_i(x, y, \theta) = 2\cos(\theta)(\cos(\theta) - \cos(\varphi)) + \sin^2(\theta) + \sin(\theta)g_i(x, y)(\cos(\theta) - 2\cos(\varphi)), \quad (3.30)$$

$$\text{for all } x \in \mathbb{R}, \ y \in (\alpha_i(x), \beta_i(x)), \ \theta \in (-\varphi, \varphi), \ i = 1,...,n$$

$$a_i(x_i, y_i, \theta_i) = \cos(\theta_i)\left(\frac{y_i - \alpha_i(x_i)}{\beta_i(x_i) - \alpha_i(x_i)}\beta_i''(x_i) + \left(1 - \frac{y_i - \alpha_i(x_i)}{\beta_i(x_i) - \alpha_i(x_i)}\right)\alpha_i''(x_i)\right)$$
$$+ \left((\beta_i(x_i) - \alpha_i(x_i))\sin(\theta_i) - g_i(x_i, y_i)\cos(\theta_i)\right)\frac{\beta_i'(x_i) - \alpha_i'(x_i)}{(\beta_i(x_i) - \alpha_i(x_i))^2} \quad (3.31)$$

for all $x \in \mathbb{R}$, $y \in (\alpha_i(x), \beta_i(x))$, $\theta \in (-\varphi, \varphi)$, $i = 1, \ldots, n$.

Condition (3.8) guarantees that $h_i(x, y, \theta) > 0$ for all $x \in \mathbb{R}$, $y \in (\alpha_i(x), \beta_i(x))$, $\theta \in (-\varphi, \varphi)$ $i = 1, \ldots, n$ with $\varphi \in \left(0, \frac{\pi}{2}\right)$ that satisfies (3.4), see Lemma 1 in Section 4.

**Discussion of the cruise controller (3.26), (3.27):** (i) The cruise controller (3.26), (3.27) for each vehicle requires measurement of its own state $(x_i, y_i, \theta_i, v_i)$ and the "relative" positions and relative speeds $(x_j - x_i, y_j - y_i, \dot{x}_i - \dot{x}_j, \dot{y}_i - \dot{y}_j)$ of all vehicles that are within distance $d_{i,j} = \sqrt{(x_i - x_j)^2 + p_{i,j}(y_i - y_j)^2}$ less than $\lambda$. Moreover, the cruise controller (3.26), (3.27) for each vehicle requires knowledge of the values of its own boundary functions $\alpha_i(x)$, $\beta_i(x)$ and their derivatives $\alpha_i'(x), \beta_i'(x)$ and $\alpha_i''(x), \beta_i''(x)$ at the current x-coordinate of the position of the vehicle. Therefore, the cruise controller (3.26), (3.27) is completely decentralized and does not require any a priori knowledge of the corridor sets $Y_i$ defined by (3.5) (only local knowledge and information is required).

(ii) The design of the cruise controller (3.26), (3.27) is based on the Lyapunov-like function $H$ defined by (3.21). This function has many differences with the Lyapunov functions that were used in [15], [16] since here we consider the general case where each vehicle has its own desired speed and its own corridor.

The following result guarantees that the system is well-posed, namely the solution is defined for all $t \geq 0$.

**Theorem 2:** *For every $w_0 \in \Omega$, the initial value problem (3.3), (3.26), (3.27) with initial condition $w(0) = w_0$ has a unique solution $w(t)$, defined for all $t \geq 0$, that satisfies $w(t) \in \Omega$ for all $t \geq 0$.*

**Discussion of Theorem 2:** Theorem 2 guarantees that (i) no collisions will occur, (ii) no vehicle will exit the corridor set $Y_i = \{(x, y) : x \in \mathbb{R}, \alpha_i(x) < y < \beta_i(x)\}$ in which it is moving, (iii) no vehicle will have a speed greater than the speed limit, (iv) no vehicle will move in the upstream direction, and (v) no vehicle will have orientation greater than $\varphi$. However, Theorem 2 does not guarantee a positive lower bound for the speed of all vehicles, i.e., it may happen that $\lim_{t \to +\infty}(v_i(t)) = 0$ for some $i \in \{1, \ldots, n\}$, and the corresponding vehicles will tend to stop. On the other hand, Theorem 3 below guarantees that this is not the case when the "boundary" functions $\alpha_i(x), \beta_i(x)$ are constant and the speed set-points are equal. Simulations show (see next section) that the case where some of the vehicles tend to stop may occur only for very special initial conditions and for very special geometries.

The proof of Theorem 2 is a consequence of Theorem 1 and the following lemma.

**Lemma 2**: *Let constants $v_{max} > 0$, $v_i^* \in (0, v_{max})$, $L_{i,j} > 0$, $i, j = 1,...,n$, $i \neq j$, $\lambda > 0$ that satisfies (3.16), and $\varphi \in \left(0, \frac{\pi}{2}\right)$ that satisfies (3.4) be given. Define the function $H : \Omega \to \mathbb{R}_+$ by means of (3.21), where $\Omega$ is given by (3.12). Moreover, for every $r > 0$ define the set*

$$S_r := \{w \in \Omega : H(w) \leq r\} \tag{3.32}$$

*Then, there exist constants $A > 1$, $\bar{\varphi} \in (0, \varphi)$, $\underline{v} \in \left(0, \min_{i=1,...,n}\{v_i^*\}\right)$, $\bar{v} \in \left(\max_{i=1,...,n}\{v_i^*\}, v_{max}\right)$, $\delta \in (0,1)$ and $\rho > 0$, such that for all $w \in S_r$ and $i, j = 1,...,n$, with $j \neq i$, the following inequalities hold*

$$d_{i,j} \geq AL_{i,j} \tag{3.33}$$

$$\frac{1}{2}\left((1+\delta)\alpha_i(x_i) + (1-\delta)\beta_i(x_i)\right) \leq y_i \leq \frac{1}{2}\left((1-\delta)\alpha_i(x_i) + (1+\delta)\beta_i(x_i)\right) \tag{3.34}$$

$$|\theta_i| \leq \bar{\varphi} \tag{3.35}$$

$$\underline{v} \leq v_i \leq \bar{v} \tag{3.36}$$

$$\begin{aligned}&|\Phi_i(w)| + |\Xi_i(w)| + |F_i| + |u_i| + \left|\frac{d}{dt}(\Phi_i(w))\right| + \left|\frac{d}{dt}(\Xi_i(w))\right| \\ &+ \frac{1}{q_i(v_i, \theta_i, f_i(w))} + \sum_{j \neq i}\left|\frac{d}{dt}(d_{i,j})\right| + \left|\frac{d}{dt}\left(U'\left(\frac{2y_i - (\beta_i(x_i) + \alpha_i(x_i))}{\beta_i(x_i) - \alpha_i(x_i)}\right)\right)\right| \leq \rho\end{aligned} \tag{3.37}$$

Indeed, if we denote by $f(w)$ the vector $\dot{w}$ for the closed-loop system (3.3), (3.26), (3.27), then Lemma 2 and the state space $\Omega \subset \mathbb{R}^{4n}$ guarantee that the set $f(S_r)$ is bounded for every $r > 0$. Although the function $H(w)$ cannot bound a size function for the closed-loop system (3.3), (3.26), (3.27), we can apply Theorem 1 and guarantee forward completeness.

Our second result deals with the case where paths of all vehicles are in the same set defined by $\alpha_i(x) \equiv \zeta_i$, $\beta_i(x) \equiv \eta_i$ for $i = 1,...,n$, and all vehicles have the same speed set-points $v_i^* = v^*$ for $i = 1,...,n$. In this case, the solution is not only well-defined, but, in addition, the speed $v_i$ of each vehicle converges (asymptotically) to the speed set-point $v^*$ and the orientation $\theta_i$ and accelerations of each vehicle converge (asymptotically) to 0.

**Theorem 3:** *Suppose that $\beta_i(x) \equiv \eta_i \in \mathbb{R}$ and $\alpha_i(x) \equiv \zeta_i \in \mathbb{R}$ with $r_{max} \geq \eta_i - \zeta_i \geq r_{min} > 0$ and $v_i^* = v^*$ for all $i = 1,...,n$. Then, for every solution $w(t)$ of the closed-loop system (3.3), (3.26), (3.27), the following equations hold for $i = 1,...,n$:*

$$\lim_{t \to +\infty}(v_i(t)) = v^*, \quad \lim_{t \to +\infty}(\theta_i(t)) = 0 \tag{3.38}$$

$$\lim_{t \to +\infty}(\Phi_i(w(t))) = \lim_{t \to +\infty}\left(U_i'\left(\frac{y_i(t) - (\eta_i + \zeta_i)}{\eta_i - \zeta_i}\right)\frac{2}{\eta_i - \zeta_i} + \Xi_i(w(t))\right) = 0 \tag{3.39}$$

$$\lim_{t \to +\infty}(F_i(t)) = 0, \quad \lim_{t \to +\infty}(u_i(t)) = 0 \tag{3.40}$$

**Discussion of Theorem 3: (i)** The proof of Theorem 3 is based on Barbalat's Lemma, and thus only asymptotic convergence is guaranteed.

**(ii)** When $\beta_i(x) \equiv \eta_i$, $a_i(x) \equiv \zeta_i$, and $v_i^* = v^* \in (0, v_{\max})$ for all $i = 1,...,n$, the closed-loop system (3.3), (3.26), (3.27) under the change of coordinates $\tilde{x}_i(t) = x_i(t) - \int_0^t \left( n^{-1} \sum_{j=1}^n v_j(s) \cos(\theta_j(s)) \right) ds - n^{-1} \sum_{j=1}^n x_j(0)$ is given by (after a slight abuse of notation – dropping the tildes),

$$\dot{x}_i = v_i \cos(\theta_i) - n^{-1} \sum_{j=1}^n v_j \cos(\theta_j)$$
$$\dot{y}_i = v_i \sin(\theta_i)$$
$$\dot{\theta}_i = u_i$$
$$\dot{v}_i = F_i$$
(3.41)

System (3.41) evolves on the set $\Omega$ defined by (3.12) and satisfies $\dfrac{d}{dt}\left(\sum_{i=1}^n x_i\right) = 0$. Thus, the state of (3.41) satisfies the additional condition $\sum_{i=1}^n x_i = 0$ (center of mass transformation; we have brought the x-coordinate of the center of mass of the vehicles to be at 0). The set of equilibrium points of the transformed closed-loop system (3.41), (3.26), (3.27) is described by

$$E = \left\{ w \in \Omega, \begin{array}{l} v_i - v^* = \theta_i = \Phi_i(w) = 0, i = 1,...,n \\ U_i'\left(\dfrac{2y_i - (\eta_i + \zeta_i)}{(\eta_i - \zeta_i)}\right) \dfrac{2}{(\eta_i - \zeta_i)} + \Xi_i(w) = 0, i = 1,...,n \\ \sum_{i=1}^n x_i = 0 \end{array} \right\}$$

Notice that if there exists a constant $M > 0$ for which the solution of the closed-loop system (3.3), (3.26), (3.27) satisfies $|x_i(t) - x_j(t)| \leq M$ for all $i, j = 1,...,n$ and $t \geq 0$, then the corresponding solution of the closed-loop system (3.41), (3.26), (3.27) is bounded. If we denote by $f(w)$ the vector $\dot{w}$ defined by (3.41), (3.26) and (3.27) then we can state that Theorem 2 guarantees that for every solution $w(t)$ of the closed-loop system (3.41), (3.26), (3.27), it holds that $\lim_{t \to +\infty}(\dot{w}(t)) = \lim_{t \to +\infty}(f(w(t))) = 0$. If for some initial condition $w_0 \in \Omega$ with $\sum_{i=1}^n x_i(0) = 0$, the corresponding solution $w(t)$ is bounded, then there exists a compact set $K \subset \Omega$ such that $w(t) \in K$ for all $t \geq 0$. Consequently, in this case the omega limit set $\omega(w_0)$ (see [20]) is non-empty and satisfies $\omega(w_0) \subseteq E$. Since $\lim_{t \to +\infty}(dist(w(t), \omega(w_0))) = 0$ (see [20]) and $\omega(w_0) \subseteq E$ (which implies that $dist(w(t), \omega(w_0)) \leq dist(w(t), E)$ we can guarantee in the case that $\lim_{t \to +\infty}(dist(w(t), E)) = 0$. Thus, the controllers (3.26), (3.27) guarantee that every bounded solution of the closed-loop system (3.41), (3.26), (3.27) tends to the set $E$ as $t \to +\infty$ (see Theorem 3 above). However, it should be noted that the limit $\lim_{t \to +\infty}(dist(w(t), E)) = 0$ does not imply the existence of a point $w^* \in E$ for which $\lim_{t \to +\infty}(|w(t) - w^*|) = 0$, i.e., we cannot guarantee that an "ultimate" arrangement exists for the vehicles. On the other hand, simulations (see next section) show that in all cases that we have tested, the solutions of the closed-loop system (3.41), (3.26), (3.27) are bounded and there is a point $w^* \in E$ (that depends on the initial condition) for which $\lim_{t \to +\infty}(\|w(t) - w^*\|) = 0$, i.e., the vehicles tend to a final arrangement.

(iii) Theorem 3 guarantees that if there exists a constant $r \in \mathbb{R}$ for which $\beta_i(x) \equiv \eta_i$ and $\alpha_i(x) \equiv \zeta_i$ for all $x \geq r \in \mathbb{R}$, namely, if the boundaries of the corridor of each vehicle remain constant for $x \geq r$ and if in addition there exists a time $T > 0$ for which $\min_{i=1,...,n}(x_i(T)) \geq r$, then the speed of each vehicle converges to the common speed set point $v^*$ (see (3.38)) and the orientation of each vehicle converges to 0 (see (3.40)). This case is illustrated in the simulation examples of the following section.

## 4. Numerical Simulations

To illustrate the effectiveness of the controllers (3.26), (3.27), we consider three different scenarios on lane-free roads, [26]. The first scenario demonstrates that under very special initial conditions and specific geometry, the speed of the vehicles tends to zero (see Discussion on Theorem 2). The second scenario demonstrates that the cruise controllers can be used in presence of off-ramps by appropriate selection of the boundary functions $\alpha_i$ and $\beta_i$ of exiting vehicles, while the third scenario considers the case of road narrowing and convergence of the speed of all vehicles to a common speed set-point.

In all cases, the vehicle-repulsive potential function $V$ and the boundary-repulsive potential function $U$ are defined by

$$V_i(d) = \begin{cases} q \dfrac{(\lambda - d)^3}{d - L}, & L < d \leq \lambda \\ 0 & d > \lambda \end{cases}, \quad i = 1,...,n \quad (4.1)$$

$$U_i(y) = \begin{cases} \left(\dfrac{1}{1-y^2} - c\right)^2, & y \in \left(-1, -\dfrac{\sqrt{c-1}}{c}\right) \cup \left(\dfrac{\sqrt{c-1}}{c}, 1\right) \\ 0 & , y \in \left[-\dfrac{\sqrt{c-1}}{c}, \dfrac{\sqrt{c-1}}{c}\right] \end{cases}, \quad i = 1,...,n \quad (4.2)$$

where $c \geq 1, q > 0$ are design parameters. Both $V$ and $U$ above, satisfy (3.13) - (3.15) and (3.17), (3.17), respectively. The constant $q$, can be used to adjust the repulsion force of the potential $V$ in (4.1) and, consequently, the magnitude of the acceleration $F_i$, see (3.26). The constant $c \geq 1$ affects the boundary repulsion and condition (3.17). More specifically, for $c > 1$, we have that $U(y) = 0$ in a neighborhood around $y = 0$. In the third scenario (constant boundary functions $\alpha_i(x) = \eta$, $\beta_i(x) = \zeta$ and $v_i^* = v^*$, $i = 1,...,n$), the constant $c$ may affect also the ultimate arrangement of the vehicles. In the numerical simulations below, the following parameters were used (unless otherwise specified) $\mu_1 = 0.1$, $\mu_2 = 1$, $R = 1$, $c = 2.1$, $q = 10^{-5}$, $\varepsilon = 0.001$, $L = 6$, $p_{i,j} = 4.25$, $i,j = 1,...,n$, $\varphi = 0.45$, $v_{\max} = 35$, $\lambda = 100$.

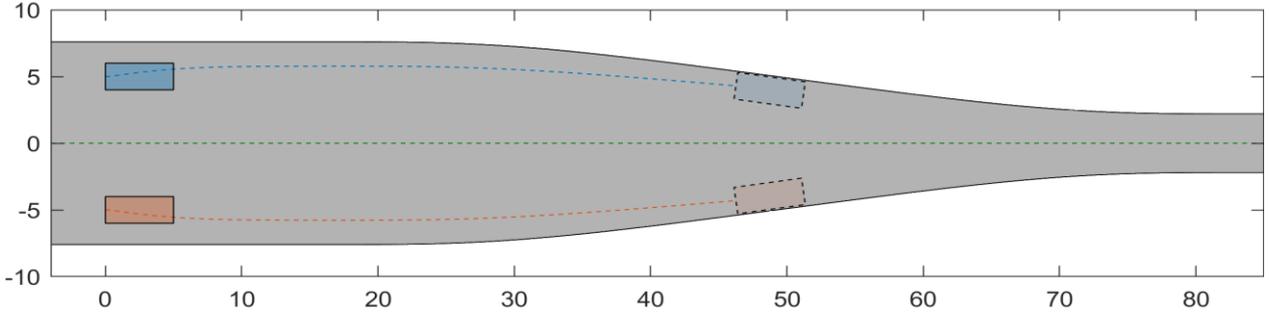
**Figure 1:** Initial position of vehicles and road configuration.

### 4.1. *Topological Obstructions.*

We consider two vehicles with $x_i(0)=0$, $v_i(0)=30$, $\theta_i(0)=0$, $i=1,2$, and $y_1(0)=5$, $y_2(0)=-5$ in a road section, as shown in Figure 1, that satisfies $\beta(x)=-\alpha(x)$ and $\beta(x)-\alpha(x)=4$, for $x>80$. Moreover, we select $p_{i,j}=1$, $i,j=1,2$, $L=6$, and $v^*=30$. Notice that both the road boundaries and the position of the vehicles are completely symmetric with respect to $y=0$, and for $x>80$ the width of the road is small enough so that only one vehicle can be placed laterally. Figure 2, shows the speeds of the vehicles which both converge to zero. This scenario illustrates that in very specific cases, Theorem 2 does not guarantee a lower bound on the speeds. Note that by slightly changing the initial position of the vehicles $x_1(0)=0$ and $x_2(0)=0.2$, both vehicles can pass through the bottleneck, see Figure 3, while their speeds converge to $v^*$ as shown in Figure 4 (the road has constant width for $x>80$ - recall Theorem 3). Videos of both situations may be watched for better appreciation https://youtu.be/aQ2KMcg7UMo

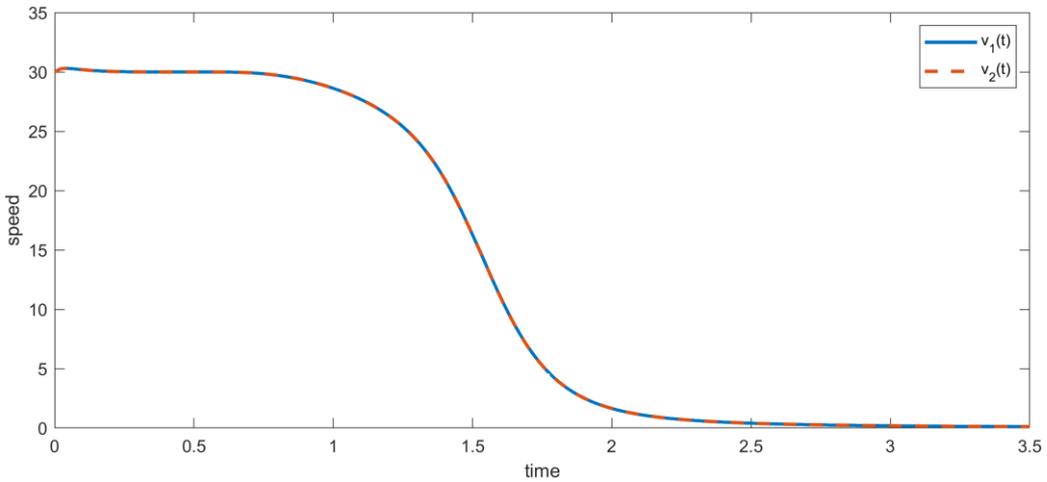
**Figure 2:** Speed of vehicles converging to zero

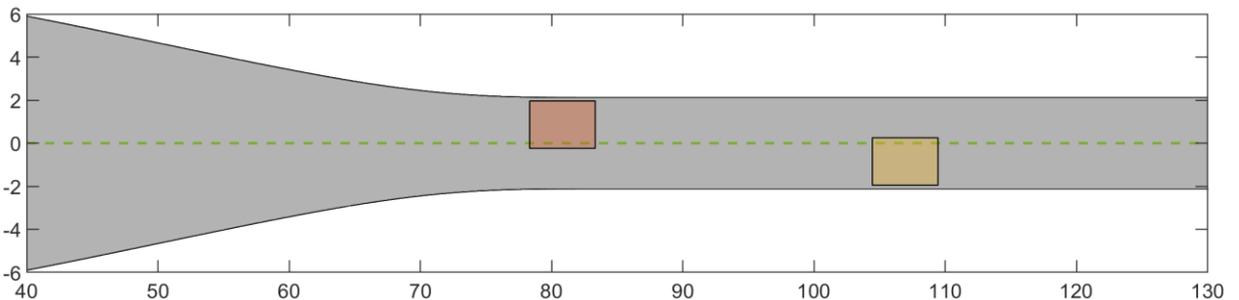
**Figure 3:** Vehicle position at $t=3.2$

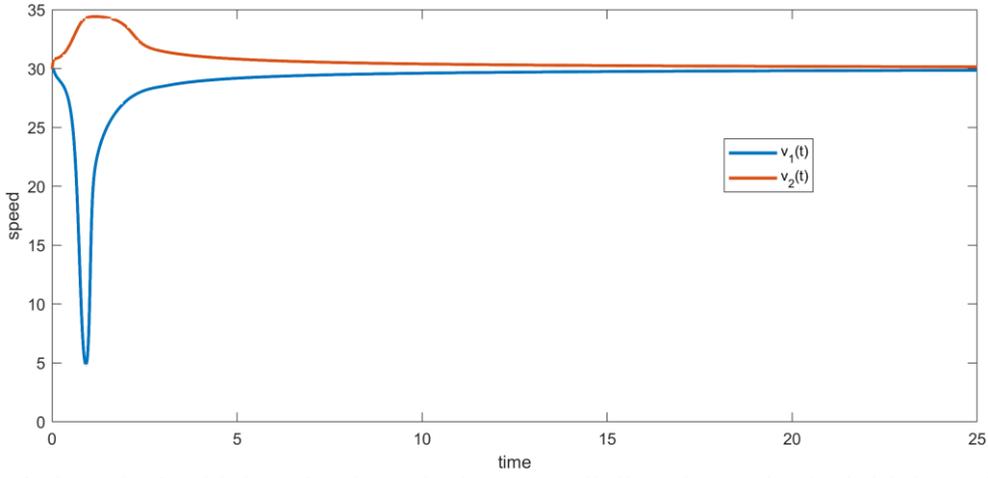

**Figure 4:** Speed of vehicles after introducing a small disturbance in the initial conditions.

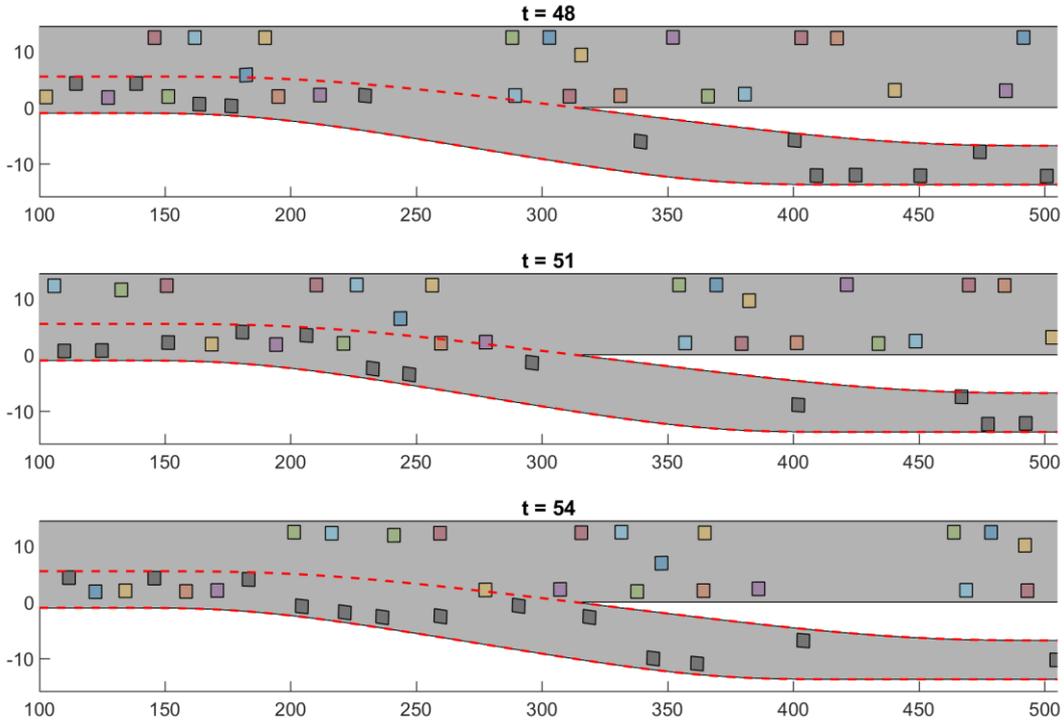

**Figure 5:** Simulation layout for various time instances.

### 4.2. Road with Off-ramp

In this scenario we consider $n=150$ vehicles on a lane-free road with an off-ramp. The corridor boundaries $\alpha_i(x)$ and $\beta_i(x)$ were appropriately selected for certain vehicles to exit the highway. Vehicles start randomly scattered on the road with zero orientation and equal speeds of $v_i(0)=30$. Moreover, vehicles exiting the highway are marked with black colour and follow the corridor marked with red dashed lines towards the off-ramp. All other vehicles follow the main road. Finally, the speed set-point $v_i^*$ for each vehicle is selected randomly from the set $\{28,29,30,31\}$. It should be noted that we have not observed any deadlocks as in the special case of Section 4.1. Figure 5 shows the simulation layout at the off-ramp at different time instances. The minimum inter-vehicle distance is shown in Figure 6, illustrating that there are no collisions between vehicles since $\min_{i,j=1,\ldots,n}(d_{i,j}(t))>L$. Vehicles exiting at the off-ramp separate swiftly from mainstream vehicles, as can also be watched in the corresponding video https://youtu.be/h4nIUz2nl9s.

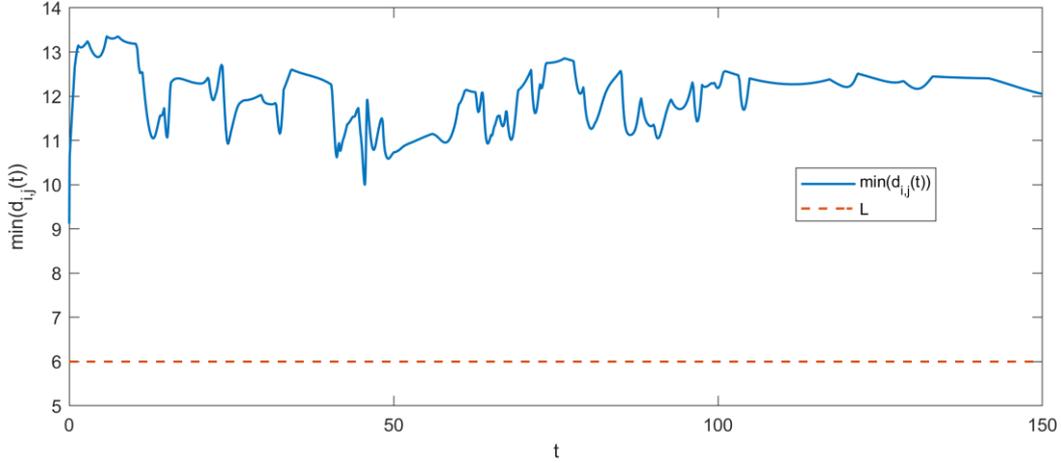

**Figure 6:** Minimum inter-vehicle distance

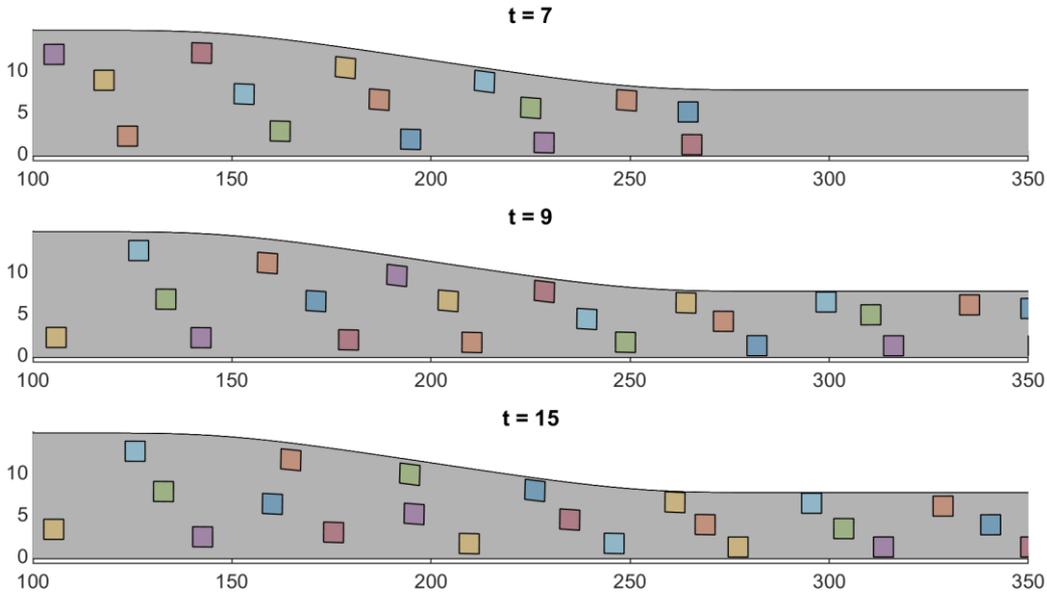

**Figure 7:** Simulation layout of a narrowing road.

### 4.3. Road Narrowing

We consider $n=50$ vehicles on a lane-free road whose initial width of 14.4 reduces by 50% and remains the same thereafter, namely $\alpha_i(x) \equiv 0$ and $\beta_i(x) \equiv 7.2$ for $x \geq 300$, $i=1,...,30$. Moreover, the speed set-point is the same for all vehicles, $v_i^* = v^* = 30$ for $i=1,...,n$. Figure 7 shows the layout of the narrowing of the road at various times where all vehicles swiftly approach and squeeze into the bottleneck while avoiding collision with each other and the boundaries of the road. Since the width of the road remains constant for $x \geq 300$, the speed of all vehicles converges to the common speed set-point $v^*$, (see Discussion of Theorem 3). Figure 8 shows the convergence of $|(v_1(t)-v^*,...,v_{30}(t)-v^*)|_\infty$ confirming that the speeds of the vehicles converge to the speed set-point (conclusion of Theorem 3). Finally, Figure 9 shows the minimum inter-vehicle distance. On one hand, it is shown that there are no collisions among vehicles ($\min_{i=1,...,n}(d_{i,j}(t)) > L$), and on the other hand that the inter-vehicle distance converges to a constant value, indicating that the vehicles tend to reach an ultimate arrangement. It should be noted again, that no vehicles enter a deadlock where their speeds converge to zero. Again, vehicles swiftly approach and squeeze into the bottleneck before accelerating and re-establishing longer distances, as can also be watched in the corresponding video https://youtu.be/1b3hp8v9juw.

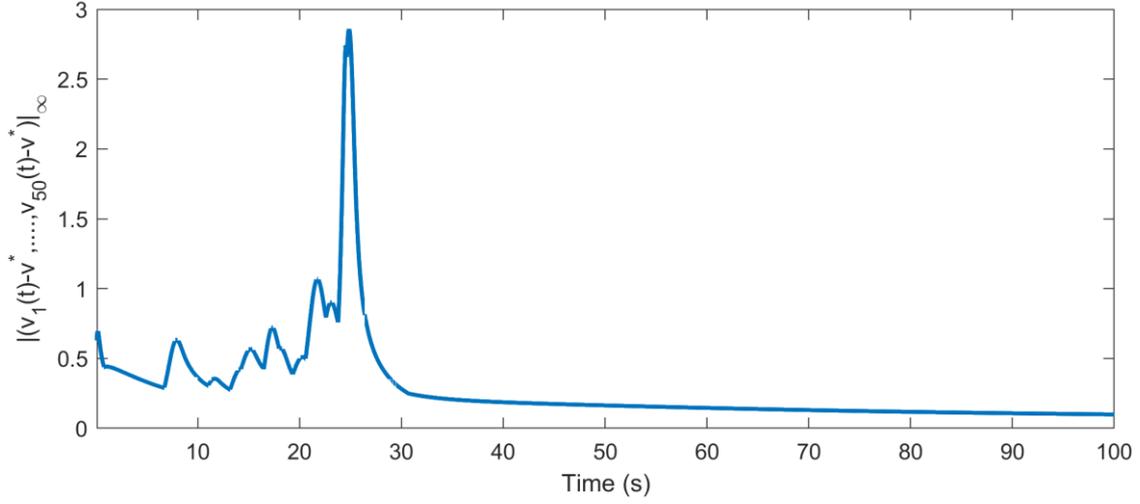

**Figure 8:** Convergence of $|(v_1(t)-v^*,...,v_{50}(t)-v^*)|_\infty$

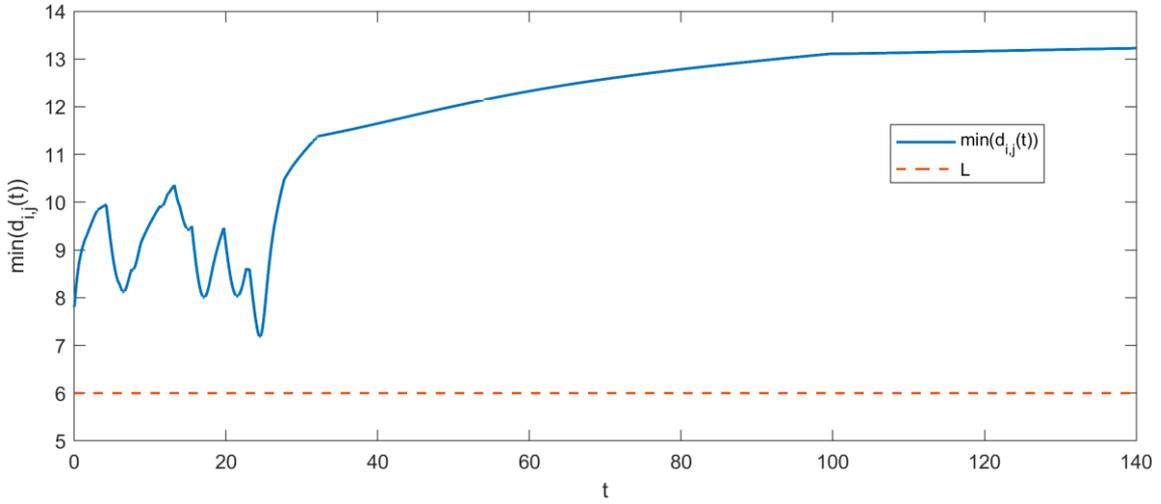

**Figure 9:** Minimum inter-vehicle distance.

## 5. Proofs

For the proof of Theorem 2 we require the following technical lemma.

**Lemma 1:** *Let* $\varphi \in \left(0, \frac{\pi}{2}\right)$ *such that* $\cos(\varphi) > \frac{1}{3}$, *and* $\alpha_i(x)$, $\beta_i(x)$ *satisfying (3.7) and (3.8).*
*Define* $\Gamma := \max\left(\sup_{x \in \mathbb{R}}(|\alpha_i'(x)|), \sup_{x \in \mathbb{R}}(|\beta_i'(x)|)\right)$. *Then, there exists a constant* $\varepsilon > 0$, *such that*
$|g_i(x,y)| \leq \Gamma$ *and* $h_i(x,y,\theta) \geq \varepsilon$, *for all* $x \in \mathbb{R}$, $y \in (\alpha_i(x), \beta_i(x))$, $\theta \in (-\varphi, \varphi)$, *and* $i = 1,...,n$.

**Proof of Lemma 1:** Notice first that from definition (3.22) we have

$$g_i(x,y) = \frac{y - \alpha_i(x)}{\beta_i(x) - \alpha_i(x)} \beta_i'(x) + \left(1 - \frac{y - \alpha_i(x)}{\beta_i(x) - \alpha_i(x)}\right) \alpha_i'(x) \qquad (5.1)$$

Moreover, due to (3.7) it follows that $\frac{y-\alpha_i(x)}{\beta_i(x)-\alpha_i(x)} \in (0,1)$ for all $x \in \mathbb{R}$ and $y \in (\alpha_i(x), \beta_i(x))$. The latter, together with (5.1) and the fact that $\Gamma := \max\left(\sup_{x \in \mathbb{R}}(|\alpha_i'(x)|), \sup_{x \in \mathbb{R}}(|\beta_i'(x)|)\right)$ implies that $|g_i(x,y)| \leq \Gamma$. Next, notice that inequalities $\cos(\varphi) > \frac{1}{3}$ and $1 \geq \cos(\theta) > \cos(\varphi)$, imply that

$$|2\cos(\varphi) - \cos(\theta)| < \cos(\varphi) \tag{5.2}$$

Then, from definition (3.30) and (5.2) we get

$$\begin{aligned} h_i(x,y,\theta) &= 1 + (\cos(\theta) - 2\cos(\varphi))(\cos(\theta) + g_i(x,y)\sin(\theta)) \\ &\geq 1 - \cos(\varphi)|\cos(\theta) + g_i(x,y)\sin(\theta)| \end{aligned} \tag{5.3}$$

Since $|g_i(x,y)| \leq \Gamma < \tan(\varphi)$ for all $x \in \mathbb{R}$, and $y \in (\alpha_i(x), \beta_i(x))$, $i = 1, \ldots, n$, there exists $\omega \in [-\arctan(\Gamma), \arctan(\Gamma)]$ such that $g_i(x,y) = \tan(\omega)$. Hence, using (5.3) and inequalities $\cos(\omega) \geq \cos(\arctan(\Gamma)) > \cos(\varphi)$ and $|\cos(\omega - \theta)| \leq 1$ we have that

$$h_i(x,y,\theta) \geq 1 - \cos(\varphi)|\cos(\theta) + \tan(\omega)\sin(\theta)| = 1 - \frac{\cos(\varphi)}{\cos(\omega)}|\cos(\omega - \theta)|$$

$$\geq 1 - \frac{\cos(\varphi)}{\cos(\omega)} \geq 1 - \frac{\cos(\varphi)}{\cos(\arctan(\Gamma))}$$

Therefore, the estimate $h_i(x,y,\theta) \geq \varepsilon$ holds for all $x \in \mathbb{R}$, $y \in (\alpha_i(x), \beta_i(x))$, $\theta \in (-\varphi, \varphi)$, and $i = 1, \ldots, n$ with $\varepsilon := 1 - \frac{\cos(\varphi)}{\cos(\arctan(\Gamma))} > 0$. This completes the proof. $\square$

**Proof of Lemma 2**: Let $r > 0$ be given and consider the set $S_r$ defined by (3.32). Since $V_{i,j}(d_{i,j}) \leq H(w)$ for all $w \in \Omega$ (recall (3.21)), Lemma 1 in [16] guarantees the existence of non-increasing functions $\rho_{i,j} : \mathbb{R}_+ \to (L_{i,j}, \lambda]$, $i, j = 1, \ldots, n$, $j \neq i$ with $\rho_{i,j}(s) \equiv \rho_{j,i}(s)$, such that the following implications hold:

$$d_{i,j} \geq \rho_{i,j}(H(w)), \text{ for all } w \in \Omega \text{ and } i, j = 1, \ldots, n, \ j \neq i \tag{5.4}$$

Due to (3.32) and (5.4) it holds that

$$d_{i,j} \geq \rho_{i,j}(r) \text{ for all } w \in S_r \text{ and } i, j = 1, \ldots, n, \ j \neq i. \tag{5.5}$$

Let $A = \min_{i,j=1,\ldots,n, j \neq i} \left\{ \frac{\rho_{i,j}(r)}{L_{i,j}} \right\} > 1$ (recall that $\rho_{i,j}(s) > L_{i,j}$ for all $s \geq 0$, $i, j = 1, \ldots, n$, $j \neq i$), then (3.33) is a direct consequence of (5.5) and the definition of $A$.

Define

$$B_{i,j} := \max\{|V_{i,j}'(d)| : AL_{i,j} \leq d \leq \lambda\} \text{ for } i, j = 1, \ldots, n, \ j \neq i \tag{5.6}$$

Notice that definition (3.10) implies that

$$\left|\frac{x_i - x_j}{d_{i,j}}\right| \leq 1 \text{ and } \left|\frac{y_i - y_j}{d_{i,j}}\right| \leq \frac{1}{\sqrt{p_{i,j}}}, \text{ for all } w \in \Omega \text{ and } i, j = 1, \ldots, n, \ j \neq i. \tag{5.7}$$

Moreover, for each $i=1,...,n$, let $m_i \geq 2$ be the maximum number of points that can be placed within the area bounded by two concentric ellipses with semi-major axes $\bar{L}_i = \min\{L_{i,j}, j=1,...,n, j \neq i\}$ and $\lambda$ satisfying (3.16), and semi-minor axes $\dfrac{\bar{L}_i}{\max\limits_{j \neq i} \sqrt{p_{i,j}}}$ and $\dfrac{\lambda}{\min\limits_{j \neq i} \sqrt{p_{i,j}}}$ so that each point has distance (in the metric given by (3.5)) at least $\bar{L}_i$ from every other point. Then, it follows from (3.10), (3.14), the fact that $d_{i,j} > L_{i,j}$ for $i,j=1,...,n$, $j \neq i$ and the definition of $m_i$ above, that the sums $\sum\limits_{j \neq i} V'_{i,j}(d_{i,j}) \dfrac{(x_i - x_j)}{d_{i,j}}$, $\sum\limits_{j \neq i} p_{i,j} V'_{i,j}(d_{i,j}) \dfrac{(y_i - y_j)}{d_{i,j}}$ contain at most $m_i$ non-zero terms, namely the terms with $d_{i,j} \leq \lambda$. Definition (5.6) in conjunction with (3.24), (3.25), (5.4), (5.7) and the fact that $p_{i,j} > 0$ for all $i,j=1,...,n$, $j \neq i$, gives the following estimate for all $w \in \Omega$ and $i=1,...,n$:

$$\max\{|\Phi_i(w)|, |\Xi_i(w)|\} \leq \max\left\{\sum_{j \neq i} |V'_{i,j}(d_{i,j})| \left|\dfrac{x_i - x_j}{d_{i,j}}\right|, \sum_{j \neq i} p_{i,j} |V'_{i,j}(d_{i,j})| \left|\dfrac{y_i - y_j}{d_{i,j}}\right|\right\}$$

$$\leq \max\left\{\sum_{j \neq i} |V'_{i,j}(d_{i,j})|, \sum_{j \neq i} \sqrt{p_{i,j}} |V'_{i,j}(d_{i,j})|\right\} \leq \sum_{j \neq i} \sqrt{(1+p_{i,j})} B_{i,j} \leq \max_{j \neq i}\{m_i \sqrt{(1+p_{i,j})} B_{i,j}\} \quad (5.8)$$

$$\leq \max_{\substack{i=1,...,n \\ j \neq i}}\{m_i \sqrt{(1+p_{i,j})} B_{i,j}\} =: \xi$$

We show next that there exists $\delta \in (0,1)$ such that (3.34) holds. Notice first that due to (3.17), (3.18) it follows that there exists a constant $\delta \in (0,1)$ for which $\{s \in (-1,1): U_i(s) \leq r\} \subset [-\delta, \delta]$. The latter, together with the fact that $U_i\left(\dfrac{2y_i - (\beta_i(x_i) + \alpha_i(x_i))}{\beta_i(x_i) - \alpha_i(x_i)}\right) \leq H(w) \leq r$, $i=1,...,n$ for all $w \in S_r$ (recall (3.21) and (3.32)) implies that

$$1 - \delta \leq \dfrac{2(y_i - \alpha_i(x_i))}{\beta_i(x_i) - \alpha_i(x_i)} \leq 1 + \delta$$

or equivalently

$$\dfrac{1}{2}\left((1+\delta)\alpha_i(x_i) + (1-\delta)\beta_i(x_i)\right) \leq y_i \leq \dfrac{1}{2}\left((1-\delta)\alpha_i(x_i) + (1+\delta)\beta_i(x_i)\right), \quad i=1,...,n \quad (5.9)$$

establishing inequality (3.34).

In order to show inequality (3.35), notice first that, due to definitions (3.32) and (3.21) it holds that

$$\dfrac{(\sin(\theta_i) - g_i(x_i, y_i)\cos(\theta_i))^2}{\cos(\theta_i) - \cos(\varphi)} \leq \dfrac{2r}{R} \quad (5.10)$$

Moreover, due to Lemma 1 it holds that

$$|g_i(x, y)| \leq \Gamma$$

with $\Gamma = \max\left(\sup_{x\in\mathbb{R}}(|\alpha_i'(x)|), \sup_{x\in\mathbb{R}}(|\beta_i'(x)|)\right) < \tan(\varphi)$. Therefore, there exists a constant $\varepsilon \in (0,1)$ such that $\varepsilon\tan^2(\varphi) > \Gamma^2$ or equivalently

$$\frac{1}{\varepsilon^{-1}\Gamma^2 + 1} > \cos^2(\varphi) \tag{5.11}$$

Using inequality (5.10), the fact that $|g_i(x,y)| \leq \Gamma$ and inequality $(\sin(\theta) - \lambda\cos(\theta))^2 \geq (1-\varepsilon)\sin^2(\theta) + (1-\varepsilon^{-1})\lambda^2\cos^2(\theta)$, it follows that

$$(1-\varepsilon)\sin^2(\theta_i) + (1-\varepsilon^{-1})\lambda^2\cos^2(\theta_i) \leq \kappa(\cos(\theta_i) - \cos(\varphi)) \tag{5.12}$$

for all $\lambda \in [-\Gamma, \Gamma]$, $i = 1,\dots,n$, where $\kappa = \frac{2r}{R}$. Inequality (5.12) can be further reduced to

$$(\varepsilon^{-1}\lambda^2 + 1)\cos^2(\theta_i) + \frac{\kappa}{1-\varepsilon}\cos(\theta_i) - \left(\frac{\kappa}{1-\varepsilon}\cos(\varphi) + 1\right) \geq 0 \tag{5.13}$$

Define the function

$$B(\lambda) := \frac{\sqrt{\frac{\kappa^2}{(1-\varepsilon)^2} + 4(\varepsilon^{-1}\lambda^2 + 1)\left(\frac{\kappa}{1-\varepsilon}\cos(\varphi) + 1\right)} - \frac{\kappa}{1-\varepsilon}}{2(\varepsilon^{-1}\lambda^2 + 1)} \text{, for } \lambda \in [-\Gamma, \Gamma] \tag{5.14}$$

for $\varepsilon \in (0,1)$, $\kappa = 2\frac{r}{R}$ and notice that inequality (5.13) is satisfied for all $w \in S_r$ for which it holds

$$\cos(\theta_i) \geq B(\lambda), \; i = 1,\dots,n. \tag{5.15}$$

Moreover, it should be noticed that due to (5.11) satisfies $1 > B(\lambda)$ for all $\lambda \in \mathbb{R}$ and $B(\lambda) \geq B(\Gamma) > \cos(\varphi)$ for all $\lambda \in [-\Gamma, \Gamma]$. Consequently, there exists $\bar{\varphi} \in (0, \varphi)$ such that $\cos(\bar{\varphi}) = B(\Gamma)$ which due to (5.15) implies that

$$|\theta_i| \leq \bar{\varphi} \tag{5.16}$$

establishing inequality (3.35)

Finally, we show that there exist constants $\underline{v} \in \left(0, \min_{i=1,\dots,n}\{v_i^*\}\right)$ and $\bar{v} \in \left(\max_{i=1,\dots,n}\{v_i^*\}, v_{\max}\right)$ such that (3.36) holds for all $w \in S_r$. Due to (3.32) and (3.21) it holds that

$$\frac{(v_i\cos(\theta_i) - f_i(w))^2}{v_i(v_{\max} - v_i)} \leq 2r \;,\; i = 1,\dots,n \tag{5.17}$$

for all $w \in S_r$. Inequality (5.8), the facts that $|g_i(x,y)| \leq \Gamma$, $i = 1,\dots,n$ for all $w \in \Omega$ and $b$ is non-increasing imply that $b(\Phi_i(w) + g_i(x_i, y_i)\Xi_i(w)) \geq b(\xi(1+\Gamma))$ for all $w \in S_r$. Since $b(x) \in (0,1]$ for all $x \in \mathbb{R}$, definition (3.23) gives that $f_i(w) \in [\underline{\sigma}, \bar{\sigma}]$, for all $w \in S_r$ with $\underline{\sigma} = b(\xi(1+\Gamma))\min_{i=1,\dots,n}(v_i^*)$, $\bar{\sigma} = \max_{i=1,\dots,n}(v_i^*)$. Inequalities (5.17) give for $i = 1,\dots,n$:

$$(2r + \cos^2(\theta_i))v_i^2 - 2(f_i(w)\cos(\theta_i) + rv_{\max})v_i + f_i^2(w) \leq 0 \tag{5.18}$$

Notice that since $v_{\max}\cos(\varphi) > v_i^*$ (recall (3.4)), $\cos(\theta_i) > \cos(\varphi)$ for $\theta_i \in (-\varphi,\varphi)$, and $f_i(w) \in [\underline{\sigma},\overline{\sigma}]$, we get

$$4\left(f_i(w)\cos(\theta_i) + rv_{\max}\right)^2 - 4f_i^2(w)\left(2r + \cos^2(\theta_i)\right) = 4r\left(rv_{\max}^2 + 2f_i(w)\left(v_{\max}\cos(\theta_i) - f_i(w)\right)\right)$$
$$> 4r\left(rv_{\max}^2 + 2f_i(w)\left(v_{\max}\cos(\varphi) - f_i(w)\right)\right) > 0$$

for all $\theta_i \in (-\varphi,\varphi)$. It follows then, by solving inequalities (5.18) with respect to $v_i$ that

$$\frac{f_i^2(w)}{f_i(w)\cos(\theta_i) + rv_{\max} + \sqrt{r\left(rv_{\max}^2 + 2v_{\max}f_i(w)\cos(\theta_i) - 2f_i^2(w)\right)}} \leq v_i$$

$$v_i \leq \frac{f_i(w)\cos(\theta_i) + rv_{\max} + \sqrt{r\left(rv_{\max}^2 + 2v_{\max}f_i(w)\cos(\theta_i) - 2f_i^2(w)\right)}}{2r + \cos^2(\theta_i)}$$
(5.19)

Define

$$\psi(x,y) = \frac{xy + rv_{\max} + \sqrt{r\left(rv_{\max}^2 + 2v_{\max}xy - 2x^2\right)}}{2r + y^2}, \quad x \in [\underline{\sigma},\overline{\sigma}], \ y \in [\cos(\varphi),1] \quad (5.20)$$

Notice that since $v_{\max}\cos(\varphi) > v_i^*$ for $i = 1,...,n$ we get $\overline{\sigma} = \max_{i=1,...,n}(v_i^*) < v_{\max}\cos(\varphi)$, which implies that $\psi(x,y) < v_{\max}$ for all $x \in [\underline{\sigma},\overline{\sigma}]$, $y \in [\cos(\varphi),1]$. It follows then from continuity of the function $\psi(x,y)$ on $[\underline{\sigma},\overline{\sigma}] \times [\cos(\varphi),1]$ that $\zeta = \max\{\psi(x,y) : x \in [\underline{\sigma},\overline{\sigma}], y \in [\cos(\varphi),1]\} < v_{\max}$. Consequently, we get from (5.20) for all $\theta_i \in (-\varphi,\varphi)$ and $f_i(w) \in [\underline{\sigma},\overline{\sigma}]$, $w \in S_r$:

$$0 < \frac{\underline{\sigma}}{v_{\max}\left(1 + r + \sqrt{r(r+2)}\right)} \leq v_i \leq \zeta \quad (5.21)$$

Therefore, (3.36) holds for all $w \in S_r$ with $\underline{v} = \dfrac{\underline{\sigma}}{v_{\max}\left(1 + r + \sqrt{r(r+2)}\right)}$ and $\overline{v} = \zeta$.

We finally show that (3.37) holds. First, we prove that $\dfrac{d}{dt}(\Phi_i(w))$ is bounded for all $i = 1,...,n$. To that end, we show first that $\dot{d}_{i,j}$ is bounded for all $i,j = 1,...,n$ with $j \neq i$. Indeed, using (3.3), the Cauchy-Schwarz inequality, and definition (3.10), the assumption that $p_{i,j} > 0$, for $i,j = 1,...,n$, $j \neq i$, inequality $-2v_iv_j\cos(\theta_i - \theta_j) \leq v_i^2 + v_j^2$ (which holds for all $v_i \in (0,v_{\max})$) and the fact that $v_i \in (0,v_{\max})$, we obtain the following estimate

$$|\dot{d}_{i,j}| \leq \sqrt{\left(v_i\cos(\theta_i) - v_j\cos(\theta_j)\right)^2 + p_{i,j}\left(v_i\sin(\theta_i) - v_j\sin(\theta_j)\right)^2}$$
$$\leq \sqrt{(1 + p_{i,j})}\sqrt{v_i^2 + v_j^2 - 2v_iv_j\cos(\theta_i - \theta_j)} \quad (5.22)$$
$$\leq \sqrt{(1 + p_{i,j})}(v_i + v_j) \leq 2\sqrt{2(1 + p_{i,j})}v_{\max}$$

Due to (3.3) and definition (3.24) we have

$$\frac{d}{dt}\big(\Phi_i(w)\big) = \sum_{j \neq i} V''_{i,j}(d_{i,j})\dot{d}_{i,j} \frac{(x_i - x_j)}{d_{i,j}} - \sum_{j \neq i} V'_{i,j}(d_{i,j})\dot{d}_{i,j} \frac{(x_i - x_j)}{d_{i,j}^2} + \sum_{j \neq i} V'_{i,j}(d_{i,j}) \frac{(\dot{x}_i - \dot{x}_j)}{d_{i,j}} \quad (5.23)$$

It follows from (3.14) and (3.33) that $V'_{i,j}(d_{i,j})$, $V''_{i,j}(d_{i,j})$ are bounded for all $i, j = 1, ..., n$ with $j \neq i$. The latter, in conjunction with (3.3), (3.33), (5.7), (5.16), (5.21), (5.22), and formula (5.23), implies that $\frac{d}{dt}\big(\Phi_i(w)\big)$ is bounded for all $i = 1, ..., n$. Similarly, we can prove that $\frac{d}{dt}\big(\Xi_i(w)\big)$ is bounded for all $i = 1, ..., n$.

We show now that $u_i$ is bounded for $i = 1, ..., n$. Since $\frac{y - \alpha_i(x)}{\beta_i(x) - \alpha_i(x)} \in (0,1)$ for all $x \in \mathbb{R}$ and $y \in (\alpha_i(x), \beta_i(x))$, we get by using (3.7), (3.8), (3.9), Lemma 1, definition (3.31), and inequality $|\tan(\theta)| \leq \tan(\varphi)$ for $\theta \in (-\varphi, \varphi)$ that

$$|a_i(x, y, \theta)| \leq \sup_{x \in \mathbb{R}}\big(|\alpha''_i(x)| + |\beta''_i(x)|\big) + (r_{\max} + 1)\frac{2\tan^2(\varphi)}{r_{\min}^2} \quad (5.24)$$

for all $x \in \mathbb{R}$, $y \in (\alpha_i(x), \beta_i(x))$, $\theta \in (-\varphi, \varphi)$ and $i = 1, ..., n$. Moreover, due to Lemma 1, (5.8), (5.9), (5.16), (5.21), inequalities $|\sin(\theta)| < \sin(\varphi)$, $\cos(\theta) > \cos(\varphi)$, $|\sin(\theta) - g_i(x, y)\cos(\theta)| < 2\tan(\varphi)$ (a direct consequence of Lemma 1), definition (3.27), (3.7), (5.24), and inequality

$$\left| U'_i\left(\frac{2y_i - (\beta_i(x_i) + \alpha_i(x_i))}{\beta_i(x_i) - \alpha_i(x_i)}\right) \right| \leq \max\big\{|U'_i(s)| : s \in [-\delta, \delta]\big\} := J_i \quad (5.25)$$

we obtain that

$$|u_i| \leq \frac{v_{\max}(1 - \cos(\varphi))^2}{\varepsilon R}\left(2\mu_2 v_{\max}\tan(\varphi) + R\frac{|a_i(x_i, y_i, \theta_i)|}{\cos(\overline{\varphi}) - \cos(\varphi)} + \frac{2J_i}{r_{\min}} + \xi\right) \quad (5.26)$$

implying that $u_i$ are bounded for $i = 1, ..., n$.

We show now that $F_i$ is bounded for $i = 1, ..., n$. Using (3.4), (3.23), the fact that $b$ is non-increasing with $b(x) \leq 1$ for all $x \in \mathbb{R}$, we get $\cos(\theta_i) > \cos(\varphi) > \frac{v_i^*}{v_{\max}} \geq \frac{f_i(w)}{v_{\max}}$ for all $\theta_i \in (-\varphi, \varphi)$, $i = 1, ..., n$. Thus, due to (3.23), (5.21), (5.8), the fact that $b$ is non-increasing, and inequality $|g_i(x, y)| \leq \Gamma$, it follows that

$$v_{\max} v_i \cos(\theta_i) + v_{\max} f_i(w) - 2f_i(w)v_i \geq v_i^* b(\xi(1+\Gamma))(v_{\max} - \overline{v}) > 0 \quad (5.27)$$

Inequalities (5.27) and definition (3.29), give the following estimate for $i = 1, ..., n$

$$\frac{1}{q(v_i, \theta_i, f_i(w))} \leq \frac{2v_{\max}^3}{v_i^* b(\xi(1+\Gamma))(v_{\max} - \overline{v})} =: \overline{q}_i \quad (5.28)$$

We prove that $Z_i(w)$ is bounded for $i = 1, ..., n$. We have from (3.3) and definition (3.22) that

$$\frac{d}{dt}\big(g_i(x_i,y_i)\big) = v_i \cos(\theta_i)\left(\frac{y_i - \alpha_i(x_i)}{\beta_i(x_i) - \alpha_i(x_i)}\beta_i''(x_i) + \left(1 - \frac{y_i - \alpha_i(x_i)}{\beta_i(x_i) - \alpha_i(x_i)}\right)\alpha_i''(x_i)\right)$$

$$+ v_i\big(\sin(\theta_i) - g_i(x_i,y_i)\cos(\theta_i)\big)\frac{\beta_i'(x_i) - \alpha_i'(x_i)}{\beta_i(x_i) - \alpha_i(x_i)}$$

Since $\frac{y_i - \alpha_i(x_i)}{\beta_i(x_i) - \alpha_i(x_i)} \in (0,1)$ for all $x_i \in \mathbb{R}$ and $y_i \in (\alpha_i(x_i), \beta_i(x_i))$, Lemma 1, (3.7), (3.8), (3.9), (5.21), and inequality $|\sin(\theta) - g_i(x,y)\cos(\theta)| < 2\tan(\varphi)$, we can conclude that $\frac{d}{dt}(g_i(x_i,y_i))$ is bounded for $i=1,\ldots,n$. Using (5.8), (3.28), (3.20), Lemma 1, the fact that $v_i^* \in (0, v_{\max})$ for $i=1,\ldots,n$, boundedness of $\frac{d}{dt}(g_i(x_i,y_i))$, $\frac{d}{dt}(\Phi_i(w))$, and $\frac{d}{dt}(\Xi_i(w))$, and formula

$$\frac{d}{dt}\big(\Phi_i(w) + g_i(x_i,y_i)\Xi_i(w)\big) = \frac{d}{dt}(\Phi_i(w)) + \frac{d}{dt}(\Xi_i(w))g_i(x_i,y_i) + \frac{d}{dt}(g_i(x_i,y_i))\Xi_i(w)$$

we conclude that $Z_i(w)$ is bounded for $i=1,\ldots,n$.

Using (5.8), (5.21), (5.26), (5.28), boundedness of $Z_i(w)$, the fact that $v_i^* b(\xi) \leq f_i(w) \leq v_i^*$ (recall (3.23), and that $b$ is non-increasing), and definition (3.26) we can conclude that

$$|F_i| \leq \overline{q}_i\left(\mu_1\big(v_{\max} - v_i^* b(\xi)\big) + (\tan(\varphi)+1)\xi + \frac{|Z_i(w)| + \overline{v}|u_i|}{\underline{v}(v_{\max} - \overline{v})}\right) \qquad (5.29)$$

is bounded for $i=1,\ldots,n$. Finally, using (3.3) and definition (3.22) we obtain that

$$\frac{d}{dt}\left(U_i'\left(\frac{2y_i - (\beta_i(x_i) + \alpha_i(x_i))}{\beta_i(x_i) - \alpha_i(x_i)}\right)\right) = U_i''\left(\frac{2y_i - (\beta_i(x_i) + \alpha_i(x_i))}{\beta_i(x_i) - \alpha_i(x_i)}\right)\frac{(\sin(\theta_i) - g_i(x_i,y_i)\cos(\theta_i))v_i}{\beta_i(x_i) - \alpha_i(x_i)}$$

which is bounded due to (5.21), (3.7), (3.18), (5.9), and inequality $|\sin(\theta) - g_i(x,y)\cos(\theta)| < 2\tan(\varphi)$. Inequalities (5.26), (5.28), (5.29), (5.8), (5.22), and boundedness of $\frac{d}{dt}(\Phi_i(w))$, $\frac{d}{dt}(\Xi_i(w))$, and $\frac{d}{dt}\left(U_i'\left(\frac{2y_i - (\beta_i(x_i) + \alpha_i(x_i))}{\beta_i(x_i) - \alpha_i(x_i)}\right)\right)$ for $i=1,\ldots,n$, imply that there exists a constant $\rho > 0$ that satisfies (3.37). This completes the proof. □

**Proof of Theorem 2:** Using (3.3), (3.10), (3.21) and definitions (3.30), (3.31), it follows that for all $w \in \Omega$

$$\dot{H}(w) = \sum_{i=1}^{n} \frac{(v_i \cos(\theta_i) - f_i(w))}{2(v_{max} - v_i)^2 v_i^2} (v_{max} v_i \cos(\theta_i) + f_i(w) v_{max} - 2 f_i(w) v_i) F_i$$

$$- \sum_{i=1}^{n} (v_i \cos(\theta_i) - f_i(w)) \frac{1}{v_i(v_{max} - v_i)} \frac{d}{dt}(f_i(w)) - \sum_{i=1}^{n} \frac{(v_i \cos(\theta_i) - f_i(w))}{v_i(v_{max} - v_i)} v_i \sin(\theta_i) u_i$$

$$+ \sum_{i=1}^{n} (\sin(\theta_i) - \cos(\theta_i) g_i(x_i, y_i)) \times$$

$$\left( \frac{Rh_i(x_i, y_i, \theta_i)}{2(\cos(\theta_i) - \cos(\varphi))^2} u_i - R \frac{v_i \cos(\theta_i) a_i(x_i, y_i, \theta_i)}{\cos(\theta_i) - \cos(\varphi)} + U_i' \left( \frac{2 y_i - (\beta_i(x_i) + \alpha_i(x_i))}{\beta_i(x_i) - \alpha_i(x_i)} \right) \frac{2 v_i}{(\beta_i(x_i) - \alpha_i(x_i))} \right)$$

$$+ \sum_{i=1}^{n} v_i \cos(\theta_i) \sum_{j \neq i} V'_{i,j}(d_{i,j}) \frac{(x_i - x_j)}{d_{i,j}} + \sum_{i=1}^{n} v_i \sin(\theta_i) \sum_{j \neq i} p_{i,j} V'_{i,j}(d_{i,j}) \frac{(y_i - y_j)}{d_{i,j}}$$

(5.30)

Using (5.30), definitions (3.24), (3.25), (3.28), (3.29), and adding and subtracting terms we further obtain:

$$\dot{H}(w) = \sum_{i=1}^{n} (v_i \cos(\theta_i) - f_i(w)) \left( q(v_i, \theta_i, f_i(w)) F_i - \frac{Z_i(w)}{v_i(v_{max} - v_i)} + \Phi_i(w) + g_i(x_i, y_i) \Xi_i(w) - \frac{\sin(\theta_i) u_i}{(v_{max} - v_i)} \right)$$

$$+ \sum_{i=1}^{n} (\sin(\theta_i) - \cos(\theta_i) g_i(x_i, y_i))$$

$$\times \left( \frac{Rh_i(x_i, y_i, \theta_i)}{2(\cos(\theta_i) - \cos(\varphi))^2} u_i - R \frac{v_i \cos(\theta_i) a_i(x_i, y_i, \theta_i)}{\cos(\theta_i) - \cos(\varphi)} + U_i' \left( \frac{2 y_i - (\beta_i(x_i) + \alpha_i(x_i))}{\beta_i(x_i) - \alpha_i(x_i)} \right) \frac{2 v_i}{(\beta_i(x_i) - \alpha_i(x_i))} + v_i \Xi_i(w) \right)$$

$$+ \sum_{i=1}^{n} f_i(w) (\Phi_i(w) + g_i(x_i, y_i) \Xi_i(w))$$

(5.31)

Using (3.26), (3.27), and (5.31) it follows that

$$\dot{H}(w) = -\sum_{i=1}^{n} \mu_1 (v_i \cos(\theta_i) - f_i(w))^2 - \sum_{i=1}^{n} \mu_2 v_i^2 (\sin(\theta_i) - g_i(x_i, y_i) \cos(\theta_i))^2$$

$$+ \sum_{i=1}^{n} f_i(w)(\Phi_i(w) + g_i(x_i, y_i) \Xi_i(w))$$

(5.32)

From (3.23) and (5.32), we get

$$\dot{H}(w) = -\sum_{i=1}^{n} \mu_1 (v_i \cos(\theta_i) - f_i(w))^2 - \sum_{i=1}^{n} \mu_2 v_i^2 (\sin(\theta_i) - g_i(x_i, y_i) \cos(\theta_i))^2$$

$$+ \sum_{i=1}^{n} v_i^* b (\Phi_i(w) + g_i(x_i, y_i) \Xi_i(w))(\Phi_i(w) + g_i(x_i, y_i) \Xi_i(w))$$

(5.33)

Let $z_i := \Phi_i(w) + g_i(x_i, y_i) \Xi_i(w)$, $i = 1,...,n$, and notice that due to (3.19) we have

$$\sum_{i=1}^{n} v_i^* b(z_i) z_i = \sum_{\substack{i \in \{1,\dots,n\} \\ z_i \le \varepsilon}} v_i^* b(z_i) z_i + \sum_{\substack{i \in \{1,\dots,n\} \\ z_i > \varepsilon}} v_i^* b(z_i) z_i \le v_{\max} \sum_{\substack{i \in \{1,\dots,n\} \\ z_i \le \varepsilon}} \varepsilon + v_{\max} \sum_{\substack{i \in \{1,\dots,n\} \\ z_i > \varepsilon}} M$$

$$\le v_{\max} \sum_{i=1}^{n} (\varepsilon + M) = v_{\max} n(\varepsilon + M) \tag{5.34}$$

From inequality (5.34), (5.33) and the facts that $v_i^* \in (0, v_{\max})$ and $(v_i \cos(\theta_i) - f_i(w))^2 \ge 0$, $(\sin(\theta_i) - g_i(x_i, y_i) \cos(\theta_i))^2 \ge 0$ for $i = 1,\dots,n$, $w \in \Omega$ it follows that

$$\dot{H}(w) \le n v_{\max}(M + \varepsilon) =: D \quad \text{for all } w \in \Omega \tag{5.35}$$

By virtue of Lemma 2 and Theorem 1, the solution $w(t)$ of (3.3), (3.26), (3.27) is defined for all $t \ge 0$ and satisfies $w(t) \in \Omega$ for all $t \ge 0$. The proof is complete. □

The proof of Theorem 3 is performed by using Barbălat's lemma [20] and its following variant which uses uniform continuity of the derivative of a function. For reader's convenience it is stated next.

**Lemma 3:** *If a function $\gamma \in C^2(\mathbb{R}_+)$ satisfies $\lim_{t \to +\infty}(\gamma(t)) \in \mathbb{R}$ and $\sup_{t \ge 0}(|\ddot{\gamma}(t)|) < +\infty$, then, $\lim_{t \to +\infty}(\dot{\gamma}(t)) = 0$.*

**Proof of Theorem 3:** Notice first, that by virtue of Theorem 2, every solution $w(t)$ is defined for all $t \ge 0$ and satisfies $w(t) \in \Omega$ for all $t \ge 0$. Moreover, due to (3.19), (3.23), continuity of $\Phi_i$, and the facts that $b$ is non-increasing with $b(x) \le 1$ for $x > \varepsilon$ and $\sum_{i=1}^{n} \Phi_i(w) = 0$ (a consequence of (3.15)), we get

$$\sum_{i=1}^{n} f_i(w)(\Phi_i(w)) = v^* \sum_{i=1}^{n} \Phi_i(w) b(\Phi_i(w)) = \sum_{\substack{i \in \{1,\dots,n\} \\ \Phi_i(w) \le \varepsilon}} v^* \Phi_i(w) + \sum_{\substack{i \in \{1,\dots,n\} \\ \Phi_i(w) > \varepsilon}} v^* \Phi_i(w) b(\Phi_i(w))$$

$$\le v^* \sum_{i=1}^{n} \Phi_i(w) = 0 \tag{5.36}$$

Inequalities (5.36) and (5.32) with $\beta_i(x) \equiv \eta_i$ and $\alpha_i(x) \equiv \zeta_i$ for $i = 1,\dots,n$, give that

$$\nabla H(w) \dot{w} = -\sum_{i=1}^{n} \mu_1 (v_i \cos(\theta_i) - f_i(w))^2 - \sum_{i=1}^{n} \mu_2 (v_i \sin(\theta_i))^2 + \sum_{i=1}^{n} \Phi_i(w) f_i(w) \le 0 \tag{5.37}$$

for all $w \in \Omega$, which implies that

$$H(w(t)) \le H(w_0), \text{ for all } t \in [0, +\infty) \tag{5.38}$$

We continue to prove by means of Barbalat's Lemma (see [20]) and Lemma 3 that (3.38), (3.39), and (3.40) hold. Using the following definitions

$$z_i := \frac{v_i \cos(\theta_i) - f_i(w)}{\sqrt{(v_{\max} - v_i) v_i}} \tag{5.39}$$

$$\Lambda(\theta) := \frac{\cos(\theta) - \cos(\varphi)}{R\cos(\theta)} \tag{5.40}$$

$$K(v_i) := \sqrt{(v_{\max} - v_i)v_i} \tag{5.41}$$

we obtain from (3.3), (3.26), and (3.27) that

$$\begin{aligned}\dot{z}_i &= -\mu_1 K^2(v_i) z_i - K(v_i)\Phi_i(w) \\ \dot{\theta}_i &= -v_i \Lambda(\theta_i)\left(\mu_2 v_i \sin(\theta_i) + \frac{2}{(\eta_i - \zeta_i)} U'_i\left(\frac{2y_i - (\eta_i + \zeta_i)}{\eta_i - \zeta_i}\right) + \Xi_i(w)\right)\end{aligned} \tag{5.42}$$

We show first that

$$\lim_{t\to+\infty}(\theta_i(t)) = 0, \ \lim_{t\to+\infty}(z_i(t)) = 0, \text{ for } i = 1,...,n \tag{5.43}$$

Using Lemma 2 with $r = H(w_0)$ gives the inequalities $(v_{\max} - v_i(t))v_i(t) \geq (v_{\max} - \bar{v})\underline{v}$ for $i = 1,...,n$ and due to definition (5.39), we further obtain

$$\max\left(\sum_{i=1}^n \left(\mu_1(v_{\max} - \bar{v})\underline{v} z_i^2(t)\right), \sum_{i=1}^n \mu_2 \underline{v}^2 \sin^2(\theta_i(t))\right) \leq Q(t) \tag{5.44}$$

where

$$Q(t) := \sum_{i=1}^n \mu_1\left(v_i(t)\cos(\theta_i(t)) - f_i(w(t))\right)^2 + \sum_{i=1}^n \mu_2\left(v_i(t)\sin(\theta_i(t))\right)^2 - \sum_{i=1}^n f_i(w(t))\Phi_i(w(t)) \tag{5.45}$$

Definition (5.45) and (5.37) imply that $Q(t) = -\frac{d}{dt}H(w(t)) \geq 0$ for all $t \geq 0$. Therefore, since $H(w) \geq 0$ for all $w \in \Omega$, we have:

$$\int_0^\infty Q(t)dt \leq H(w_0). \tag{5.46}$$

Since (5.46) holds, inequality (5.44) gives that

$$\max\left(\int_0^\infty z_i^2(t)dt, \int_0^\infty \sin^2(\theta_i(t))dt\right) \leq \frac{H(w_0)}{\min\left(\mu_1(v_{\max} - \bar{v})\underline{v}, \mu_2 \underline{v}^2\right)}, \text{ for } i = 1,...,n \tag{5.47}$$

Since $\dot{\theta}_i = u_i$, by Lemma 2, (5.47) and Barbalat's Lemma it holds that $\lim_{t\to+\infty}(\theta_i(t)) = 0$, $i = 1,...,n$. Thus, by virtue of (5.47) and Barbalat's Lemma, it suffices to show that $\dot{z}_i(t)$ are bounded for $i = 1,...,n$.

Notice first that definitions (5.40), and (5.41), together with Lemma 2 imply that

$$\frac{\cos(\bar{\varphi}) - \cos(\varphi)}{R} \leq \Lambda(\theta_i) \leq \frac{1 - \cos(\varphi)}{R\cos(\bar{\varphi})} \tag{5.48}$$

$$\sqrt{(v_{\max} - \bar{v})\underline{v}} \leq K(v_i) \leq \sqrt{(v_{\max} - \underline{v})\bar{v}} \tag{5.49}$$

Furthermore, from definitions (5.39), (5.41), Lemma 2, and inequality (5.49) we have that $|z_i| \leq \dfrac{v_{max} - v^* b(\rho)}{\sqrt{(v_{max} - \underline{v})\underline{v}}}$. The previous inequality, (5.42), Lemma 2, and (5.49) give that

$$|\dot{z}_i| \leq \mu_1 \frac{(v_{max} - \underline{v})\overline{v}}{\sqrt{(v_{max} - \overline{v})\underline{v}}} \left(v_{max} - v^* b(\rho)\right) + \rho\sqrt{(v_{max} - \underline{v})\overline{v}} \quad (5.50)$$

Therefore, using (5.47), (5.50) and Barbalat's Lemma, we conclude that (5.43) hold.

We next show by using Lemma 3 with $\gamma(t) = \theta_i(t)$ and $\gamma(t) = z_i(t)$, that

$$\lim_{t \to +\infty}\left(\dot{\theta}_i(t)\right) = 0, \; \lim_{t \to +\infty}\left(\dot{z}_i(t)\right) = 0 \; \text{ for } i = 1,...,n \quad (5.51)$$

Since (5.43) holds, it suffices to show that $|\ddot{\theta}_i(t)|$ and $|\ddot{z}_i(t)|$ are bounded for $i = 1,...,n$. From (5.41), Lemma 2, (5.29), and formula

$$\frac{d}{dt}(K(v_i(t))) = K'(v_i(t))\dot{v}_i(t) = \frac{(v_{max} - 2v_i(t))F_i(t)}{2\sqrt{(v_{max} - v_i(t))v_i(t)}} \quad (5.52)$$

we have that $\dfrac{d}{dt}(K(v_i(t)))$ is bounded for $i = 1,...,n$. Using (5.42), Lemma 2, (5.49), (5.52), inequality $|z_i| \leq \dfrac{v_{max} - v^* b(\rho)}{\sqrt{(v_{max} - \overline{v})\underline{v}}}$ for $i = 1,...,n$, and formula

$$\begin{aligned}\ddot{z}_i(t) = &-\mu_1 K^2(v_i(t))\dot{z}_i(t) - 2\mu_1 K(v_i(t))K'(v_i(t))F_i(t)z_i(t) \\ &-K'(v_i(t))F_i(t) - K(v_i(t))\frac{d}{dt}(\Phi_i(w(t)))\end{aligned} \quad (5.53)$$

we can conclude that $|\ddot{z}_i(t)|$ is bounded for $i = 1,...,n$.

We show next that $|\ddot{\theta}_i(t)|$ is bounded for $i = 1,...,n$. We have from definition (5.40) that

$$\frac{d}{dt}(\Lambda(\theta_i(t))) = -\frac{\cos(\varphi)\tan(\theta_i(t))}{R\cos(\theta_i(t))}\dot{\theta}_i(t) \quad (5.54)$$

which due to (3.3), Lemma 2 and inequality $|\tan(\theta)| < \tan(\varphi)$ for $\theta \in (-\varphi, \varphi)$ is bounded for $i = 1,...,n$. Combining Lemma 2, (5.42), inequality $r_{max} \geq \eta_i - \zeta_i \geq r_{min} > 0$, (3.27), (5.48), (5.54), and formula

$$\begin{aligned}\ddot{\theta}_i(t) = &-F_i(t)\Lambda(\theta_i(t))\left(\mu_2 v_i \sin(\theta_i(t)) + 2(\eta_i - \zeta_i)^{-1}U_i'\left(\frac{2y_i(t) - (\eta_i + \zeta_i)}{\eta_i - \zeta_i}\right) + \Xi_i(w(t))\right) \\ &-v_i(t)\Lambda'(\theta_i(t))\dot{\theta}_i(t)\left(\mu_2 v_i \sin(\theta_i(t)) + 2(\eta_i - \zeta_i)^{-1}U_i'\left(\frac{2y_i(t) - (\eta_i + \zeta_i)}{\eta_i - \zeta_i}\right) + \Xi_i(w(t))\right) \\ &-v_i(t)\Lambda(\theta_i(t))\left(\mu_2 F_i(t)\sin(\theta_i(t)) + \mu_2 v_i(t)\dot{\theta}_i(t)\cos(\theta_i(t))\right) \\ &-v_i(t)\Lambda(\theta_i(t))\left(2(\eta_i - \zeta_i)^{-1}\frac{d}{dt}\left(U_i'\left(\frac{2y_i(t) - (\eta_i + \zeta_i)}{\eta_i - \zeta_i}\right)\right) + \frac{d}{dt}(\Xi_i(w(t)))\right)\end{aligned} \quad (5.55)$$

we conclude that $|\ddot{\theta}_i(t)|$ is bounded for all $i=1,...,n$. Therefore, due to (5.43), (5.53), and (5.55), it follows by Lemma 3 with $\gamma(t)=z_i(t)$ and $\gamma(t)=\theta_i(t)$ that (5.51) holds.

We prove now that (3.38), (3.39), and (3.40) hold. First, notice that due to (5.42), (5.43), (5.49), and (5.51) it holds that

$$\lim_{t \to +\infty}\left(\Phi_i(w(t))\right)=0 \tag{5.56}$$

Moreover, using (5.42), (5.43), (5.48), Lemma 2, and (5.51), we also obtain that
$\lim_{t \to +\infty}\left(U_i'\left(\dfrac{2y_i(t)-(\eta_i+\zeta_i)}{\eta_i-\zeta_i}\right)\dfrac{2}{\eta_i-\zeta_i}+\Xi_i(w(t))\right)=0$ for all $i=1,...,n$, establishing (3.39).

Taking into account definitions (3.19), (3.23), and (5.56), we get that $\lim_{t \to +\infty}\left(f_i(w(t))\right)=v^*$ for $i=1,...,n$ which in conjunction with (5.39), (5.43), and Lemma 2 establishes (3.38). Notice also that due to Lemma 2, (3.20), (3.28), (5.56), and continuity of $b'(s)$ it holds that

$$\lim_{t \to +\infty}\left(Z_i(w(t))\right)=0 \tag{5.57}$$

Finally, by taking into account Lemma 2, (5.51), (5.56), (5.57), definition (3.26), and (3.38), we conclude that $\lim_{t \to +\infty}\left(F_i(t)\right)=0$, establishing (3.40). The proof is complete. □

## 6. Conclusions

We have presented a relaxed Lyapunov-like condition for forward completeness for finite dimensional systems defined on open sets that does not require boundedness of the Lyapunov-like function along the solutions of the system. Dynamical systems with constraints are typically defined on open sets, and the proposed Lyapunov-like condition provides a useful tool for the global existence of solution for all positive times. The corresponding condition is then exploited for the design of lane-free cruise controllers for the autonomous two-dimensional movement of vehicles. The derived feedback laws (cruise controllers) are decentralized and guarantee collision avoidance and boundary respect, can consider roads of variable width with on-ramps and off-ramps as well as different desired speed for each vehicle. Finally, it was shown that when all vehicles have the same speed set-point and operate on a road of constant width, their speeds converge asymptotically to the speed set-point while their acceleration and rotation rate converge (asymptotically) to zero.

While the design of the cruise controllers for roads of non-constant width guarantees that vehicles will not collide with each other, will not exit the corridor or road in which they are moving, nor have a speed greater than the speed limit, it does not guarantee a positive lower bound for the speed of all vehicles. The latter implies that in some very special cases, vehicles may asymptotically come to a halt. Simulations show that the set of initial conditions that can lead the vehicles to a halt is of zero Lebesgue measure and thus the possibility of having this outcome is zero. In future work, we will study criteria that determine whether or not the vehicles will be led to a halt.

## Acknowledgments

The research leading to these results has received funding from the European Research Council under the European Union's Horizon 2020 Research and Innovation programme/ ERC Grant Agreement n. [833915], project TrafficFluid.